\documentclass[3p,preprint]{elsarticle}

\usepackage{amssymb}  
\usepackage{amsmath}  
\usepackage{amsthm}   
\usepackage{lineno}   
\usepackage{cases}    
\usepackage{tabularx} 
\usepackage[figuresright]{rotating}  
\usepackage{color}
\usepackage{algorithm}
\usepackage{algpseudocode}
\usepackage{accents}
\usepackage{tabularx}
\usepackage{slashbox,pict2e}

\usepackage[hidelinks]{hyperref} 
\hypersetup{pdfstartview={XYZ null null 1.00}} 
\hypersetup{pdfpagemode=UseNone} 

\allowdisplaybreaks  

\journal{Elsevier}
\date{March 31, 2020}

\begin{document}

\begin{frontmatter}

\title{A Conservative High-Order Method Utilizing Dynamic Transfinite Mortar Elements for Flow Simulation on Curved Sliding Meshes}

\author{Bin Zhang\corref{cor1}}
\ead{bzh@gwmail.gwu.edu}

\author{Chunlei Liang}

\cortext[cor1]{Corresponding author.}

\address{Department of Mechanical and Aerospace Engineering, \\
  The George Washington University, Washington, DC 20052 USA}

\begin{abstract}
We present a high-order method for flow simulation on unstructured curved nonconforming sliding meshes. This method utilizes dynamic transfinite mortar elements to exchange flow information between the two sides of a sliding interface. The method is arbitrarily high-order accurate in space, provably conservative, and satisfies outflow condition. Moreover, it retains the accuracy of a time marching scheme, and thus allows substantial reduction of rotational speed effects when equipped with a high-order temporal scheme. The method's capability of simultaneously handling multiple rotational objects is also explored. Details on the implementation are provided as well.
\end{abstract}


\end{frontmatter}


\section{Introduction}
\label{sec:introduction}

Sliding mesh has been widely used in flow simulations about moving objects. For example, it is an ideal choice for handling rotational geometries such as stirred tanks \cite{bakker-1997} and helicopter rotor blades \cite{steijl-2008}. It can also be used to ensure good mesh qualities in circumstances where purely deforming mesh may otherwise be very skewed, such as for simulating oscillating wings \cite{kinsey-2008} and vortex-induced-vibration devices \cite{sarwar-2010}. In many applications, a sliding mesh method has advantages over other methods such as overset mesh methods \cite{steger-1983} and immersed boundary methods \cite{mittal-2005} for its simplicity, efficiency and accuracy. But so far, most sliding mesh methods are still limited to low-order (second order and below) schemes that are unfavorable for simulating vortex dominated flows.

Tremendous progress has been made on high-order methods in the past decades in the computational fluid dynamics community \cite{wang-2014}. For instance, such methods include the discontinuous Galerkin (DG) method \cite{reed-1973, cockburn-2011}, the spectral element method \cite{patera-1984, karniadakis-2005, kopriva-2009}, the spectral volume method \cite{wang-2002a, wang-2002b}, the spectral difference (SD) method \cite{kopriva-1996a, kopriva-1998, liu-2006, wang-2007a,balan-2012}, to name just a few. Among these methods, the SD method solves equations in differential form directly, and is one of the most efficient high-order methods. Recently, the ideas of collocating solution and flux points of the SD method and correcting fluxes using higher-degree polynomials have led to an even more efficient high-order method --- the flux reconstruction (FR) method \cite{huynh-2007, huynh-2009}, also known as the correction procedure via reconstruction (CPR) method \cite{wang-2009}. Besides its better efficiency, by choosing different correction polynomials, the FR method can recover many existing high-order schemes such as DG and SD, and can even produce new schemes that were never reported before. The stability of the FR method has been proved in \cite{jameson-2012}. New variants of the FR method have been reported in, e.g., \cite{vincent-2011, castonguay-2012, williams-2013}. The most recent developments on the FR method are summarized in \cite{wang-2016}.

With more and more applications of high-order methods to flow simulations than ever before, there is a natural need of extending the sliding mesh concept to high-order methods to tackle complex flow problems such as those mentioned above. \citet{ferrer-2012} developed a high-order sliding-mesh DG method based on modal basis functions for simulating incompressible flow problems. \citet{ramirez-2015} applied moving-least-squares stencils to the development of a high-order sliding-mesh finite volume method. The authors of the present paper developed a simple and efficient high-order sliding-mesh SD method \cite{zhang-2015b} using iso-parametric mortar elements \cite{mavriplis-1989, kopriva-1996b, kopriva-2002}, and also extended this method to sliding-deforming meshes \cite{zhang-2016a} and to handle 3D geometries \cite{zhang-2016c}. Our previous methods, however, require uniform mesh on a sliding interface, which restricts mesh generation. Recent efforts have completely lifted this restriction, and the resulting general nonuniform sliding-mesh method has been demonstrated on the FR method \cite{zhang-2018} and the SD method on hybrid grids \cite{qiu-2019}.

Our initial effort in \cite{zhang-2018} showed that transfinite mortar elements can potentially make the sliding-mesh method arbitrarily high-order accurate. In this work, we further the investigation to provide a more detailed study emphasizing the following aspects of this method: rotational speed effects on both spatial and temporal accuracies, conservation property, outflow property, free-stream preservation property, capability of dealing with multiple objects at the same time. Meanwhile, detailed steps on the implementation are also provided in this work.

The rest of this paper is organized as follows. In Section \ref{sec:equations}, we briefly describe the equations that we are going to solve numerically. Sections \ref{sec:fr} and \ref{sec:sliding} are about the numerical methods, including the FR method and the new sliding-mesh method. Verifications and applications are reported in Section \ref{sec:examples}. Finally, Section \ref{sec:conclusion} concludes this paper.

\section{The Flow Equations}
\label{sec:equations}

We numerically solve the two-dimensional Navier-Stokes equations in the following conservative form,
\begin{equation}
\frac{\partial \mathbf{Q}} {\partial t} +  \frac{\partial \mathbf{F}} {\partial x} + \frac{\partial \mathbf{G}} {\partial y} = \mathbf{0},
\label{eq:physical}
\end{equation}
where $\mathbf{Q}$ is the vector of conservative variables; $\mathbf{F}$ and $\mathbf{G}$ are the flux vectors in the $x$- and the $y$-direction, respectively. Their expressions are
\begin{align}
\mathbf{Q} &= [\rho ~ \rho u ~ \rho v ~ E]^\mathsf{T},                                       \label{eq:Q} \\[0.5em]
\mathbf{F} &= \mathbf{F}_{\text{inv}}(\mathbf{Q}) + \mathbf{F}_{\text{vis}}(\mathbf{Q},\nabla \mathbf{Q}), \label{eq:F} \\[0.5em]
\mathbf{G} &= \mathbf{G}_{\text{inv}}(\mathbf{Q}) + \mathbf{G}_{\text{vis}}(\mathbf{Q},\nabla \mathbf{Q}), \label{eq:G}
\end{align}
where $\rho$ is fluid density, $u$ and $v$ are the Cartesian velocity components, and $E$ is the total energy per volume defined as
\begin{equation}
E = \frac{p}{\gamma-1} + \frac{1}{2}\rho(u^2+v^2),
\end{equation}
where $p$ is pressure, and $\gamma$ is the ratio of specific heats and is set to 1.4 in this work for ideal gas. The inviscid and the viscous flux vectors are
\begin{equation}
\mathbf{F}_{\text{inv}} = \left[
\begin{array}{c}
\rho u       \\
\rho u^2 + p \\
\rho uv      \\
(E+p)u
\end {array}
\right], \quad
\mathbf{G}_{\text{inv}} = \left[
\begin{array}{c}
\rho v       \\
\rho uv      \\
\rho v^2 + p \\
(E+p)v
\end {array}
\right],
\label{eq:FGinv}
\end{equation}
\begin{equation}
\mathbf{F}_{\text{vis}} = -\left[
\begin{array}{c}
0         \\
\tau_{xx} \\
\tau_{yx} \\
u\tau_{xx}+v\tau_{yx}+ \kappa T_x
\end {array}
\right], \quad
\mathbf{G}_{\text{vis}} = -\left[
\begin{array}{c}
0         \\
\tau_{xy} \\
\tau_{yy} \\
u\tau_{xy}+v\tau_{yy}+ \kappa T_y
\end {array}
\right],
\label{eq:FGvis}
\end{equation}
where
\begin{equation}
\tau_{ij} = \mu (u_{i,j}+u_{j,i}) + \lambda\delta_{ij}u_{k,k}, \nonumber
\end{equation}
is the shear stress tensor, $\mu$ is the dynamic viscosity, $\lambda=-2/3\mu$ based on the Stokes' hypothesis, $\delta_{ij}$ is the Kronecker delta, $\kappa$ is the thermal conductivity, and $T$ is temperature which is related to density and pressure through the ideal gas law
\begin{equation}
p=\rho \mathcal{R} T,
\end{equation}
where $\mathcal{R}$ is the gas constant.

To deal with moving grid, we employ an arbitrary-Lagrangian-Eulerian (ALE) approach, and map a moving physical domain to a fixed computational domain. Let $(t,x,y)$ denote the physical time and space, and $(\tau,\xi,\eta)$ the computational ones. Further assume that the mapping is: $t=\tau$, $x=x(\tau,\xi,\eta)$, and $y=y(\tau,\xi,\eta)$. It can be shown that the flow equations will take the following conservative form in the computational space
\begin{equation}
\frac{\partial \widetilde{\mathbf{Q}}} {\partial t} +  \frac{\partial \widetilde{\mathbf{F}}} {\partial \xi} +
\frac{\partial \widetilde{\mathbf{G}}} {\partial \eta} = \mathbf{0},
\label{eq:computational}
\end{equation}
where $\widetilde{\mathbf{Q}}$, $\widetilde{\mathbf{F}}$, and $\widetilde{\mathbf{G}}$ are the computational solution vector and flux vectors, and they are related to the physical ones as
\begin{align}
\widetilde{\mathbf{Q}} &= |\mathcal{J}| \mathbf{Q}, \label{eq:compQ} \\
\widetilde{\mathbf{F}} &= (-x_t y_\eta + y_t x_\eta) \mathbf{Q} + y_\eta \mathbf{F} - x_\eta \mathbf{G},  \label{eq:compF} \\
\widetilde{\mathbf{G}} &= (\hphantom{-} x_t y_\xi - y_t x_\xi) \mathbf{Q} - y_\xi \mathbf{F} + x_\xi \mathbf{G}. \label{eq:compG}
\end{align}
Alternatively, let $Q$, $F$, $\widetilde{F}$, etc., each denote a component (at the same position) of the corresponding boldface vector, and then the relations in (\ref{eq:compQ})-(\ref{eq:compG}) can be written in the following matrix form
\begin{equation}
\renewcommand*{\arraystretch}{1.3}
\left[\begin{array}{c} \widetilde{Q} \\ \widetilde{F} \\ \widetilde{G} \end{array} \right] = |\mathcal{J}|\mathcal{J}^{-1}
\left[\begin{array}{c} Q \\ F \\ G \end{array} \right], \label{eq:compQFG}
\end{equation}
where $\mathcal{J}$ represents the Jacobian matrix, $|\mathcal{J}|$ is its determinant, and $\mathcal{J}^{-1}$ is the inverse Jacobian matrix. These metric terms have the following expressions,
\begin{gather}
\mathcal{J} = \frac{\partial(t,x,y)}{\partial(\tau,\xi,\eta)} =
\begin{bmatrix}
1      & 0       & 0 \\
x_\tau & x_{\xi} & x_{\eta} \\
y_\tau & y_{\xi} & y_{\eta}
\end{bmatrix}, \qquad
|\mathcal{J}| = x_\xi y_\eta - x_\eta y_\xi  \label{eq:jac1},  \\[1em]
\mathcal{J}^{-1} = \frac{\partial(\tau,\xi,\eta)}{\partial(t,x,y)} =
\begin{bmatrix}
1      & 0      & 0 \\
\xi_t  & \xi_x  & \xi_y \\
\eta_t & \eta_x & \eta_y
\end{bmatrix} = \frac{1}{|\mathcal{J}|}
\begin{bmatrix}
\phantom{-} |\mathcal{J}|        & \phantom{-}0      & \phantom{-}0 \\
-x_t y_\eta + y_t x_\eta         & \phantom{-}y_\eta & -x_\eta \\
\phantom{-}x_t y_\xi - y_t x_\xi & -y_\xi            & \phantom{-} x_\xi
\end{bmatrix}. \label{eq:jac2}
\end{gather}

\section{Flux Reconstruction Method on Conforming Moving Grid}
\label{sec:fr}

\subsection{Grid Mapping}
\label{subsec:ale}

\begin{figure}[!htb]
\centering
\includegraphics[scale=1.0]{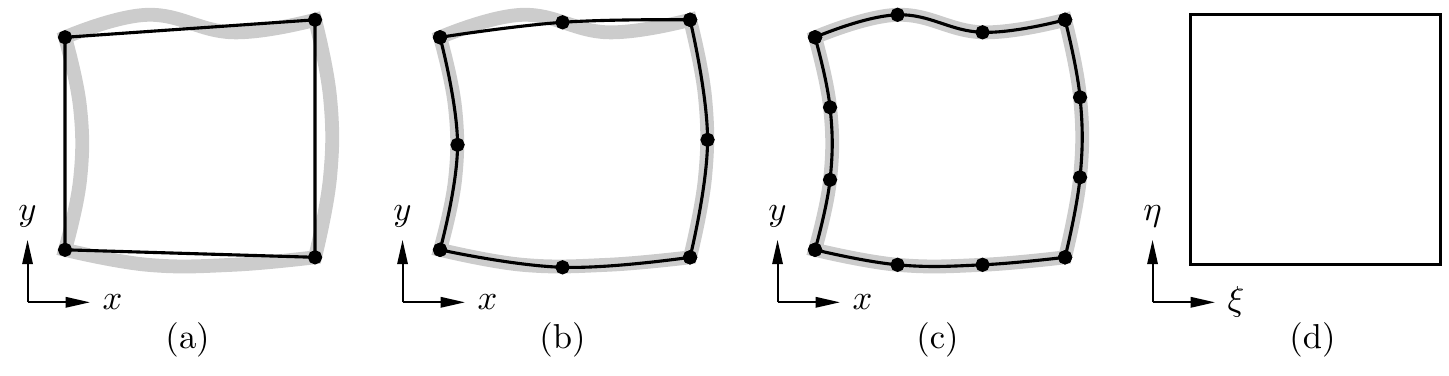}
\caption{Iso-parametric representation (black lines) of a curved element (gray lines): (a) $K=4$ (linear), (b) $K=8$ (quadratic), (c) $K=12$ (cubic), (d) a computational element (unit square).}
\label{fig:isomap2d}
\end{figure}

The first step of the FR method is to map each physical grid element to a standard computational element. In the present implementation, we discretize a computational domain into non-overlapping quadrilateral grids, and employ the following iso-parametric mapping \cite{ergatoudis-1968} to map each grid element to a unit square element (i.e., $0\leq \xi,~ \eta \leq 1$) in the computational space,
\begin{equation}
\renewcommand*{\arraystretch}{1.3}
\mathbf{x}(t,\xi,\eta) = \left[ \begin{matrix} x(t,\xi,\eta) \\ y(t,\xi,\eta) \end{matrix} \right] = \sum_{i=1}^{K} M_i(\xi,\eta) \left[ \begin{matrix} x_i(t) \\ y_i(t) \end{matrix} \right], \label{eq:isomap2d}
\end{equation}
where $K$ is the total number of nodes used to approximate a physical element, $M_i$ and $(x_i,y_i)$ are the shape function and the coordinates of the $i$-th node. Fig. \ref{fig:isomap2d} is a schematic of the iso-parametric representations of a physical element using different number of nodes. Higher-order elements (larger $K$'s) obviously represent curved boundaries more accurately. For this reason, high-order elements should be used along curved boundaries for better accuracy.

\subsection{Construction of Solution and Flux Polynomials}
\label{subsec:solpoly}

\begin{figure}[!htb]
\centering
\includegraphics[scale=1.0]{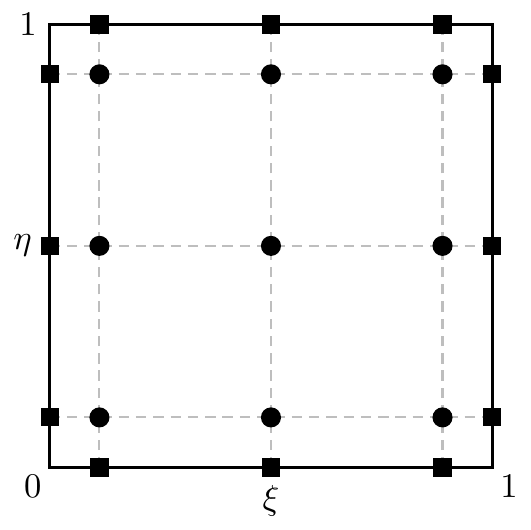}
\caption{Schematic of the distribution of SPs (circular dots) and FPs (square dots) for a third-order ($P=2$) FR scheme.}
\label{fig:spfp}
\end{figure}

Within each computational element, solution points (SPs) and flux points (FPs) are defined. The SPs are distributed along each coordinate direction inside the element, and the FPs are distributed along the boundaries only. Figure \ref{fig:spfp} is a schematic of the distribution of these points for a third-order FR scheme. In this work, the $N$ SPs/FPs in each direction of an $N$-th order scheme are chosen as the roots of the $N$-th Legendre polynomial (namely, $N$ Legendre points). At the SPs, Lagrange interpolation bases are defined. For example, the basis at the $i$-th SP along the $\xi$ direction is
\begin{align}
h_{i}(\xi) = \prod_{s=1,s\neq i}^{N}\left(\frac{\xi-X_s}{X_i-X_s}\right) \label{eq:hbasis},
\end{align}
where $X_i$ and $X_s$ are the $\xi$-coordinates of the $i$-th and the $s$-th SP, respectively. If we denote the space of all polynomials of degrees less than or equal to $N$ as $\boldsymbol{\mathsf{P}}_N$, then $h_i \in \boldsymbol{\mathsf{P}}_{N-1}$. Moreover, the $h_i$'s are linearly independent and form a basis for $\boldsymbol{\mathsf{P}}_{N-1}$.

For a function $\phi$ defined within a mesh element, assume it is smooth (e.g., $\phi \in \boldsymbol{\mathsf{C}}^\infty$) and has discrete values: $\phi_{ij}$ at $(x_i,y_j)$ that corresponds to $(X_i,X_j)$ in the computational space, where $i,j=1,2,\cdots,N$. This function can then be approximated by the following tensor-product polynomial in $\boldsymbol{\mathsf{P}}_{N-1,N-1}=\boldsymbol{\mathsf{P}}_{N-1} \otimes \boldsymbol{\mathsf{P}}_{N-1}$,
\begin{equation}
\phi(\xi,\eta) = \sum_{j=1}^{N} \sum_{i=1}^{N} \phi_{ij} h_i(\xi) h_j(\eta) \label{eq:phi},
\end{equation}
where the truncation error for the approximation is $\mathcal{O}(\xi^N,\eta^N)$ based on Taylor series expansion.

Applying (\ref{eq:phi}) to the solution and fluxes, we can obtain the following polynomial representations,
\begin{align}
\widetilde{\mathbf{Q}}(\xi,\eta) &= \sum_{j=1}^{N} \sum_{i=1}^{N} \widetilde{\mathbf{Q}}_{ij} h_i(\xi) h_j(\eta) \label{eq:Qt}, \\
\widetilde{\mathbf{F}}(\xi,\eta) &= \sum_{j=1}^{N} \sum_{i=1}^{N} \widetilde{\mathbf{F}}_{ij} h_i(\xi) h_j(\eta) \label{eq:Ft}, \\
\widetilde{\mathbf{G}}(\xi,\eta) &= \sum_{j=1}^{N} \sum_{i=1}^{N} \widetilde{\mathbf{G}}_{ij} h_i(\xi) h_j(\eta) \label{eq:Gt},
\end{align}
where the $()_{ij}$'s are the discrete values at the $(i,j)$-th SP in a standard element.

The polynomials in (\ref{eq:Qt})-(\ref{eq:Gt}) are continuous within each element, but discontinuous across element interfaces (or boundaries). For this reason, common values need to be defined at the interfaces. For instance, the common solution is computed as
\begin{equation}
\mathbf{Q}^{\text{com}} = \frac{1}{2}(\mathbf{Q}^\text{L} + \mathbf{Q}^\text{R}), \label{eq:Qcom}
\end{equation}
where $\mathbf{Q}^\text{L}$ and $\mathbf{Q}^\text{R}$ are the solution vectors on the left side and the right side of an interface.

To compute the common (normal) inviscid flux, we employ a Riemann solver, such as the following Rusanov solver \cite{rusanov-1961} with modification for moving grid,
\begin{equation}
\mathbf{F}_{\text{inv}}^{\text{com}} = \frac{1}{2} \big[(\accentset{\leftrightarrow}{\mathbf{F}}_{\text{inv}}^{\text{L}} + \accentset{\leftrightarrow}{\mathbf{F}}_{\text{inv}}^{\text{R}}) \cdot \mathbf{n} - \lambda(\mathbf{Q}^\text{R} - \mathbf{Q}^\text{L}) \big] - (\mathbf{v}_g\cdot \mathbf{n}) \mathbf{Q}^{\text{com}}, \label{eq:Fcom}
\end{equation}
where $\accentset{\leftrightarrow}{\mathbf{F}}_{\text{inv}}^{\text{L}} = (\mathbf{F}_{\text{inv}}^{\text{L}}, \mathbf{G}_{\text{inv}}^{\text{L}})$ and $\accentset{\leftrightarrow}{\mathbf{F}}_{\text{inv}}^{\text{R}} = (\mathbf{F}_{\text{inv}}^{\text{R}}, \mathbf{G}_{\text{inv}}^{\text{R}})$ are the inviscid flux vectors on the two sides of an interface; $\mathbf{n} = \mathbf{N}/\Vert \mathbf{N} \Vert$ is the unit normal vector with
\begin{equation}
\mathbf{N} = (y_\eta,-x_\eta) \text{ or } (-y_\xi, x_\xi)
\end{equation}
depending on which direction (i.e., $\xi$ or $\eta$) the interface is mapped to; $\lambda$ is the largest characteristic speed with the following expression,
\begin{equation}
\lambda= |V_n|+c = \Big|\big(\frac{1}{2}(\mathbf{v}^\text{L}+\mathbf{v}^\text{R}) - \mathbf{v}_g \big) \cdot \mathbf{n} \Big|+c,
\end{equation}
where $\mathbf{v}^\text{L}=(u^\text{L},v^\text{L})$ and $\mathbf{v}^\text{R}=(u^\text{R},v^\text{R})$ are flow velocities on the two sides of an interface; $\mathbf{v}_g=(x_t,y_t)$ represents grid velocity; $c$ is the local speed of sound. The physical normal flux in (\ref{eq:Fcom}) is converted to a computational one by multiplying the magnitude of normal, i.e.,
\begin{equation}
\widetilde{\mathbf{F}}_{\text{inv}}^{\text{com}} \text{ or } \widetilde{\mathbf{G}}_{\text{inv}}^{\text{com}} = \Vert \mathbf{N} \Vert \mathbf{F}_{\text{inv}}^{\text{com}}. \label{eq:Fn2Ft}
\end{equation}

The common viscous flux is calculated from the common solution and common gradient,
\begin{equation}
\mathbf{F}_{\text{vis}}^{\text{com}} = \mathbf{F}_{\text{vis}}(\mathbf{Q}^{\text{com}}, (\nabla\mathbf{Q})^{\text{com}}), \label{eq:Fvcom}
\end{equation}
where the common gradient is the average of the left and the right values, i.e.,
\begin{equation}
(\nabla \mathbf{Q})^{\text{com}} = ((\nabla \mathbf{Q})^\text{L} + (\nabla \mathbf{Q})^\text{R})/2.
\label{eq:nablaq_com}
\end{equation}
The procedure for calculating $\nabla \mathbf{Q}$ will be briefly discussed in the next section. The common viscous flux in (\ref{eq:Fvcom}) is converted to a computational one in the same way as (\ref{eq:Fn2Ft}).

\subsection{Flux Reconstruction}
The spatial derivatives in (\ref{eq:computational}) reduces the two flux terms to $\boldsymbol{\mathsf{P}}_{N-2,N-1}$ and $\boldsymbol{\mathsf{P}}_{N-1,N-2}$, making them inconsistent with the solution term which is $\boldsymbol{\mathsf{P}}_{N-1,N-1}$. To overcome this issue, the flux polynomials need to be reconstructed to be at least $\boldsymbol{\mathsf{P}}_{N,N-1}$ and $\boldsymbol{\mathsf{P}}_{N-1,N}$.

To do this, we use higher degree correction functions/polynomials. The corrected/reconstructed fluxes take the following form,
\begin{alignat}{5}
\widehat{\mathbf{F}}(\xi,\eta) &= \widetilde{\mathbf{F}}(\xi,\eta) &&+ \big[\widetilde{\mathbf{F}}^\text{com}(0,\eta) - \widetilde{\mathbf{F}}(0,\eta)\big] &&\cdot g_\text{\tiny L}(\xi) &&+ [\widetilde{\mathbf{F}}^\text{com}(1,\eta) - \widetilde{\mathbf{F}}(1,\eta)] &&\cdot g_\text{\tiny R}(\xi), \label{eq:Fc} \\[0.5em]
\widehat{\mathbf{G}}(\xi,\eta) &= \widetilde{\mathbf{G}}(\xi,\eta) &&+ [\widetilde{\mathbf{G}}^\text{com}(\xi,0) - \widetilde{\mathbf{G}}(\xi,0)] &&\cdot g_\text{\tiny L}(\eta) &&+ [\widetilde{\mathbf{G}}^\text{com}(\xi,1) - \widetilde{\mathbf{G}}(\xi,1)] &&\cdot g_\text{\tiny R}(\eta), \label{eq:Gc}
\end{alignat}
where $\widetilde{\mathbf{F}}(\xi,\eta)$ and $\widetilde{\mathbf{G}}(\xi,\eta)$ are the original flux polynomials from (\ref{eq:Ft}) and (\ref{eq:Gt}); $g_\text{\tiny L}$ and $g_\text{\tiny R}$ are the left and the right correction functions with degrees no less than $N$, and they are required to at least satisfy
\begin{equation}
\begin{alignedat}{2}
g_\text{\tiny L}(0) &= 1, \quad g_\text{\tiny L}(1) &&= 0, \\
g_\text{\tiny R}(0) &= 0, \quad g_\text{\tiny R}(1) &&= 1,
\end{alignedat}
\end{equation}
which ensures that
\begin{alignat}{2}
\widehat{\mathbf{F}}(0,\eta) &= \widetilde{\mathbf{F}}^\text{com}(0,\eta), \quad \widehat{\mathbf{F}}(1,\eta) &&= \widetilde{\mathbf{F}}^\text{com}(1,\eta), \\[0.5em]
\widehat{\mathbf{G}}(\xi,0) &= \widetilde{\mathbf{G}}^\text{com}(\xi,0), \quad \widehat{\mathbf{G}}(\xi,1) &&= \widetilde{\mathbf{G}}^\text{com}(\xi,1),
\end{alignat}
i.e., the reconstructed fluxes still take the common values at cell interfaces. Huynh \cite{huynh-2007} proposed several correction functions, and we employ the $g_{\text{\tiny DG}}$ correction function in this study.

Note that to calculate solution gradients consistently, the correction procedure is also applied to the solution polynomial along each direction (only for calculating the gradients). In this way, the resulting gradient polynomials are in $\boldsymbol{\mathsf{P}}_{N-1,N-1}$ as well.

\subsection{Time Marching}
With proper boundary conditions applied, the discretized equations can be written in the following residual form,
\begin{equation}
\left. \frac{\partial \widetilde{\mathbf{Q}}} {\partial t} \right|_{ij} = - \left[ \frac{\partial \widehat{\mathbf{F}}} {\partial \xi} +
\frac{\partial \widehat{\mathbf{G}}} {\partial \eta} \right]_{ij} = \mathfrak{R}_{ij}, \quad i,j=1,2,\cdots,N, \label{eq:resid}
\end{equation}
where $\mathfrak{R}_{ij}$ is the residual at the $(i,j)$-th SP. This system can be time marched using either explicit or implicit schemes. In this work, we employ several of the explicit strong stability preserving (SSP) Runge-Kutta schemes reported in \cite{kraaijevanger-1991,ruuth-2006,gottlieb-2006} for the purpose. Furthermore, in this work, all boundary conditions are weakly imposed to increase stability, and more detail can be found in, e.g., \cite{kopriva-1996a}.

\subsection{Free-Stream Preservation}
Ideally, a moving grid should not disturb a flow field. The simplest situation is that a free-stream flow must stay constant all the time on a moving grid. This is called free-stream preservation. By substituting a constant flow solution into the flow equations in (\ref{eq:computational}), we can get a system of equations that are purely about the geometrics, and are thus known as the geometric conservation law (GCL) \cite{thomas-1979},
\begin{numcases}{}
\frac{\partial (|\mathcal{J}|\xi_x)}{\partial\xi} + \frac{\partial (|\mathcal{J}|\eta_x)}{\partial\eta} = 0, \label{eq:gcl1}   \\[0.3em]
\frac{\partial (|\mathcal{J}|\xi_y)}{\partial\xi} + \frac{\partial (|\mathcal{J}|\eta_y)}{\partial\eta} = 0, \label{eq:gc2}    \\[0.3em]
\frac{\partial |\mathcal{J}|}{\partial t} + \frac{\partial (|\mathcal{J}|\xi_t)}{\partial\xi} + \frac{\partial (|\mathcal{J}|\eta_t)}{\partial\eta} = 0. \label{eq:gcl3}
\end{numcases}
To numerically satisfy free-stream preservation, the same numerical schemes for discretizing the flow equations must be applied to the GCL equations. Since the spatial discretization operator in the FR method is direct differentiation (which is exact), and the geometric terms are from analytical mapping, the first two GCL equations are hence automatically satisfied. But the third GCL equation generally can not be satisfied automatically. This problem come from the temporal discretization, for instance, a multiple-stage Runge-Kutta scheme, which is not exact. To overcome this issue, we treat $|\mathcal{J}|$ as an unknown, and solve the third equation numerically. The numerical $|\mathcal{J}|$ is then used to update the physical solution according to (\ref{eq:compQ}). In this way, the GCL is numerically satisfied, and free-stream preservation is ensured. Similar approach was reported in, e.g., \cite{minoli-2011}.

\section{A Sliding-Mesh Method using Transfinite Mortar Elements}
\label{sec:sliding}
In this section, we first introduce the concepts of sliding mesh and mortar element. Following that, a brief description is given on transfinite mapping. The projection procedures between cell faces and mortars are described subsequently. After that, we prove that the present method is globally conservative and satisfies outflow condition. We also provide details on the implementation at the end of this section.

\subsection{Sliding Mesh and Mortar Elements}
\begin{figure}[!htb]
\centering
\includegraphics[scale=1.0]{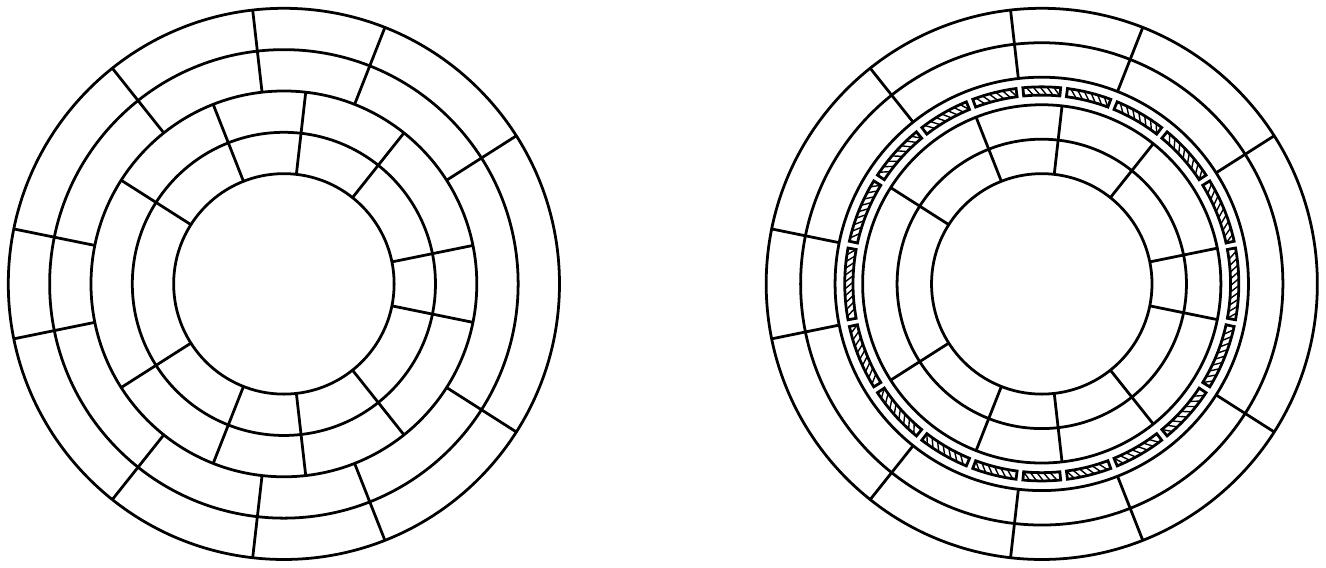}
\caption{Left: a nonconforming sliding mesh; right: distribution of mortar elements (represented by hatched lines).}
\label{fig:mortar_schematic}
\end{figure}

The simplest sliding mesh involves only two non-overlapping subdomains. Such an example is shown in Fig. \ref{fig:mortar_schematic}(a), where the two subdomains are allowed to have relative rotational motions, which leads to a nonconforming sliding interface in between. To ensure continuity of solution and conservation of fluxes on this nonconforming mesh, we employ mortar elements as the communicator between the two subdomains.

The distribution of mortar elements for this mesh is shown in Fig. \ref{fig:mortar_schematic}(b). A mortar element is formed between two successive points on the sliding interface. At every time instant, a mortar element is connected to a cell face on its left and a cell face on its right (from a counterclockwise perspective). A cell face is connected to one or multiple mortars on its one side. We call this information as the cell face and mortar connectivities. Note that these connectivities are time-dependent, and need to be updated at every stage of a time marching scheme. An efficient procedure for updating these connectivities will be presented in a later section.

\subsection{Transfinite Mapping}
\label{sec:transfinite}
\begin{figure}[!htb]
\centering
\includegraphics[scale=0.9]{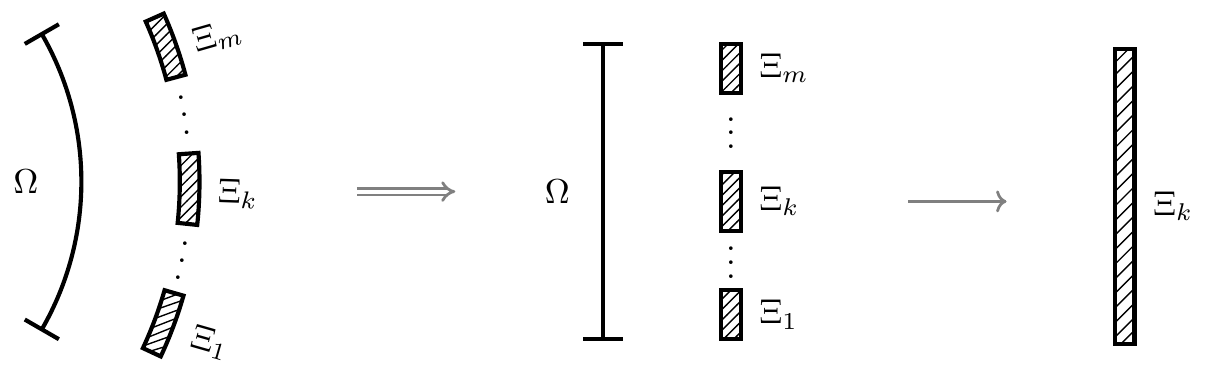}
\caption{Mapping a cell face and its mortars to unit line elements: left, physical space; middle, computational space; right, mortar space.}
\label{fig:mortar_mapping}
\end{figure}

Each sliding cell face is mapped to a unit straight line element (e.g., $0 \le \xi \le 1$) when the underlying grid cell is mapped to a unit square element in the computational space. Similarly, we also map each mortar element to a unit straight line element $0\le z \le 1$ in the mortar space to facilitate the construction of polynomials. This process has been demonstrated in Fig. \ref{fig:mortar_mapping}, where $\Omega$ denotes a cell face, and the $\Xi$'s denote the associated mortar elements. Since the shapes of a sliding cell face and a mortar element are special (i.e., they are circular arcs), we employ transfinite mapping \cite{gordon-1973a, gordon-1973b} instead of iso-parametric mapping for the task. The advantage of transfinite mapping is that it represents geometries that have analytical expressions (such as a circular arc) exactly, i.e., introduces no geometric error at all.
\begin{figure}[!htb]
\centering
\includegraphics[scale=1.0]{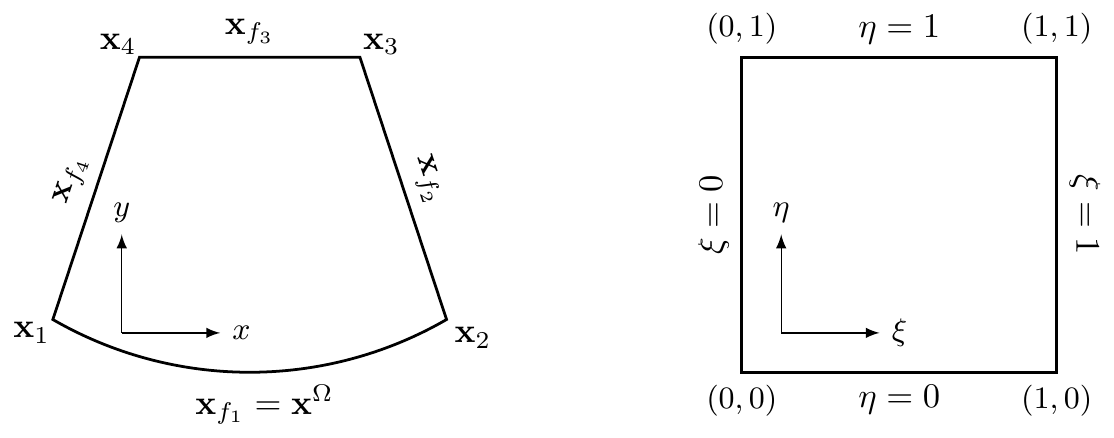}
\caption{Transfinite mapping of a physical element to a unit square element.}
\label{fig:transfinite}
\end{figure}

Figure \ref{fig:transfinite} is a schematic of a physical element and the corresponding computational element. The transfinite mapping between these two elements can be expressed as
\begin{equation}
\begin{split}
\mathbf{x}(t,\xi,\eta) = {}& (1-\eta) \mathbf{x}_{f_1}(t,\xi) + \xi \mathbf{x}_{f_2}(t,\eta) + \eta \mathbf{x}_{f_3}(t,\xi) + (1-\xi) \mathbf{x}_{f_4}(t,\eta) \\
& -(1-\xi)(1-\eta) \mathbf{x}_1(t) - \xi (1-\eta) \mathbf{x}_2(t) - \xi \eta \mathbf{x}_3(t) - (1-\xi) \eta \mathbf{x}_4(t),
\end{split}
\end{equation}
where the $\mathbf{x}_i$'s denote coordinates of the corner nodes, and $\mathbf{x}_{f_i}$'s are expressions of the faces, and these terms are all time-dependent for moving grids. If the faces are represented by one-dimensional iso-parametric mappings of the same order, then the transfinite mapping is equivalent to the iso-parametric mapping in (\ref{eq:isomap2d}). If a face has an exact expression, then that face is represented exactly by the transfinite mapping.

In our case, assume face $f_1$ is a circular arc (i.e., corresponds to the sliding cell face $\Omega$ in Fig. \ref{fig:mortar_mapping}), then it can be analytically expressed as
\begin{equation}
\renewcommand*{\arraystretch}{1.3}
\mathbf{x}^{\Omega}(t,\xi) = \left[
\begin{matrix}
x^\Omega(t,\xi) \\
y^\Omega(t,\xi)
\end{matrix} \right] =
\left[
\begin{matrix}
R \cdot \cos\big[(1-\xi)\theta^{\Omega}_1(t) + \xi \theta^{\Omega}_2(t)\big] + x_c(t) \\
R \cdot \sin\big[(1-\xi)\theta^{\Omega}_1(t) + \xi \theta^{\Omega}_2(t)\big] + y_c(t)
\end{matrix}
\right], \label{eq:face_transfinite}
\end{equation}
where $R$ and $(x_c,y_c)$ are the radius and the center coordinates of the arc; $\theta^{\Omega}_1$ and $\theta^{\Omega}_2$ are the angles corresponding to the starting point (i.e., $\mathbf{x}_1$) and the ending point (i.e., $\mathbf{x}_2$), respectively, of the face. Similarly, a mortar element $\Xi_k$ has the following exact expression,
\begin{equation}
\renewcommand*{\arraystretch}{1.3}
\mathbf{x}^{\Xi_k}(t,z) = \left[
\begin{matrix}
x^{\Xi_k}(t,z) \\
y^{\Xi_k}(t,z)
\end{matrix} \right] =
\left[
\begin{matrix}
R \cdot \cos\big[(1-z)\theta^{\Xi_k}_1(t) + z \theta^{\Xi_k}_2(t)\big] + x_c(t) \\
R \cdot \sin\big[(1-z)\theta^{\Xi_k}_1(t) + z \theta^{\Xi_k}_2(t)\big] + y_c(t)
\end{matrix}
\right], \label{eq:mortar_transfinite}
\end{equation}
where $\theta^{\Xi_k}_1$ and $\theta^{\Xi_k}_2$ are the starting and the ending angles of the mortar.

The analytical expressions in (\ref{eq:face_transfinite}) and (\ref{eq:mortar_transfinite}) are actually equivalent to the following linear mappings of the angles,
\begin{align}
\theta^{\Omega} &= (1-\xi)\theta^{\Omega}_1(t) + \xi \theta^{\Omega}_2(t), \\
\theta^{\Xi_k}    &= (1-z)\theta^{\Xi_k}_1(t) + z \theta^{\Xi_k}_2(t),
\end{align}
where $\theta^{\Omega}$ is the angle of a physical point (on $\Omega$) that is mapped to a point $\xi$ in the computational space, and $\theta^{\Xi_k}$ is the angle of a physical point (on $\Xi_k$) that is mapped to a point $z$ in the mortar space. Referring to the physical space in Fig. \ref{fig:mortar_mapping} and considering the fact that $\theta^{\Omega}=\theta^{\Xi_k}$ represents the same physical point, the following relation thus holds between the computational space and the mortar space,
\begin{equation}
\xi=o_k + s_k \cdot z, \label{eq:xi_z}
\end{equation}
where $s_k$ and $o_k$ are the scaling and the offset of the mortar $\Xi_k$ with respect to the cell face $\Omega$. More specifically,
\begin{align}
s_k &= \frac{\theta^{\Xi_k}_2 - \theta^{\Xi_k}_1}{\theta^\Omega_2 - \theta^\Omega_1} = \frac{\Delta \theta^{\Xi_k}}{\Delta \theta^\Omega} = \frac{R \cdot \Delta \theta^{\Xi_k}}{R \cdot \Delta \theta^\Omega} = \frac{L^{\Xi_k}}{L^\Omega}, \label{eq:sk} \\[1em]
o_k &= \frac{\Delta \theta^{\Xi_1} + \Delta \theta^{\Xi_2} + \cdots + \Delta \theta^{\Xi_{k-1}}}{\Delta \theta^\Omega} = \sum\nolimits_{\alpha=1}^{k-1} s_\alpha,  \label{eq:ok}
\end{align}
where $L^{\Xi_k}$ and $L^{\Omega}$ are the physical lengths of the mortar and the face. Note that the scaling and the offset are both time-dependent, and are updated at the time the connectivities are updated.

Two comparison examples of the iso-parametric mapping and the transfinite mapping are included in Appendix \hyperref[sec:appendix_transfinite]{A}.

\subsection{Projection Procedures}
Communications on a sliding interface include: projection of local discontinuous values from cell faces to mortars, computation of common values on mortars, and projection of common values back to cell faces. These procedures are discussed in details in what follows. To facilitate the discussion, we adopt the following notations: $Q$ denotes a component of $\mathbf{Q}$; $Q_i$ denotes the discrete value of $Q$ at the $i$-th FP;  $\boldsymbol{Q}$ denotes the vector $(Q_1,Q_2,\cdots,Q_N)$. This same rule also applies to fluxes.

\subsubsection{Project local values to mortars}
Take solution as an example. Each solution component on a cell face is represented by the following one-dimensional polynomial
\begin{equation}
Q^{\Omega} (\xi) = \sum_{i=1}^{N} Q_i^{\Omega}h_i(\xi). \label{eq:Qface}
\end{equation}
If we define the same set of FPs in the mortar space, then the solution polynomial on a mortar can be constructed in the same way. For example, on the left side of $\Xi_k$, the solution polynomial is
\begin{equation}
Q^{\Xi_k,\text{L}} (z) = \sum_{i=1}^{N} Q_i^{\Xi_k,\text{L}} h_i(z), \label{eq:Qmortar}
\end{equation}
where $Q_i^{\Xi_k,\text{L}}$ is a solution component at the $i$-th FP on the left side of $\Xi_k$.
\begin{figure}[!htb]
\centering
\includegraphics[scale=1.0]{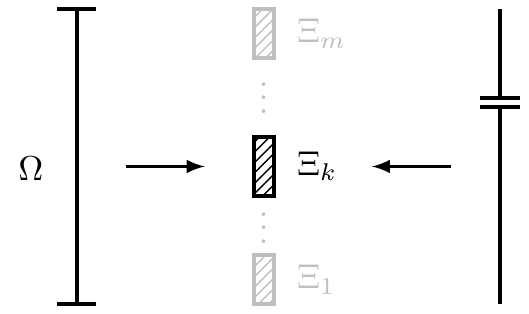}
\caption{Projection from cell faces to the two sides of a mortar.}
\label{fig:mortar_projection1}
\end{figure}
As illustrated in Fig. \ref{fig:mortar_projection1}, to get the solutions on the left side of $\Xi_k$, we require
\begin{equation}
\int_0^1 \left( Q^{\Xi_k,\text{L}}(z) - Q^{\Omega}(\xi) \right) h_j(z) ~\mathrm{d} z = 0,~~ \forall j=1,2,...,N.  \label{eq:prjQ1}
\end{equation}
Substituting (\ref{eq:Qface}) and (\ref{eq:Qmortar}) into the above equation and considering the relation in (\ref{eq:xi_z}), we will get the following equation system
\begin{equation}
\mathbf{M}\boldsymbol{Q}^{\Xi_k,\text{L}} = \mathbf{S}^{\Omega\rightarrow\Xi_k} \boldsymbol{Q}^{\Omega}, \label{eq:prjQ0}
\end{equation}
where the elements of the coefficient matrices are
\begin{align}
M_{ij} &= \int_0^1 h_i(z) h_j(z) ~\mathrm{d} z, ~~ i,j=1,2,...,N, \label{eq:mortarM} \\[0.5em]
S_{ij}^{\Omega\rightarrow\Xi_k} &= \int_0^1 h_i(o_k+s_k z) h_j(z) ~\mathrm{d} z, ~~ i,j=1,2,...,N. \label{eq:mortarS}
\end{align}
Solutions of (\ref{eq:prjQ0}), when written in matrix form, are
\begin{equation}
\boldsymbol{Q}^{\Xi_k,\text{L}} = \mathbf{P}^{\Omega\rightarrow\Xi_k} \boldsymbol{Q}^{\Omega} = \mathbf{M}^{-1} \mathbf{S}^{\Omega\rightarrow\Xi_k} \boldsymbol{Q}^{\Omega}, \label{eq:prjQ2}
\end{equation}
where $\mathbf{P}^{\Omega\rightarrow\Xi_k}$ is the projection matrix from $\Omega$ to $\Xi_k$. This process is repeated for all the solution components.

In the same way, we can get solutions on the right side of a mortar. Also note that the integrals in (\ref{eq:mortarM}) and (\ref{eq:mortarS}) can be evaluated exactly and efficiently using quadratures, e.g., the Clenshaw-Curtis quadratures. The Legendre points being the SPs/FPs in this work makes the $h_i$'s orthogonal, and in turn makes the $\mathbf{M}$ matrix diagonal and trivial to invert. The inversion actually needs to be done only once during initialization because the $\mathbf{M}$ matrix is also time-independent.

\subsubsection{Compute common values on mortars}
The common solution $\boldsymbol{Q}^{\Xi_k}$ and the inviscid normal flux $\boldsymbol{F}_{\text{inv}}^{\Xi_k}$ on a mortar are computed in the same way as on a cell interface, i.e.,
\begin{gather}
\boldsymbol{Q}^{\Xi_k} =\frac{1}{2}(\boldsymbol{Q}^{\Xi_k,\text{L}} + \boldsymbol{Q}^{\Xi_k,\text{R}}), \label{eq:Qcom_mortar} \\
\boldsymbol{F}_{\text{inv}}^{\Xi_k} = \frac{1}{2} \big[(\accentset{\leftrightarrow}{\boldsymbol{F}}_{\text{inv}}^{\Xi_k,\text{L}} + \accentset{\leftrightarrow}{\boldsymbol{F}}_{\text{inv}}^{\Xi_k,\text{R}}) \cdot \mathbf{n} - \lambda(\boldsymbol{Q}^{\Xi_k,\text{R}} - \boldsymbol{Q}^{\Xi_k,\text{L}}) \big] - (\mathbf{v}_g\cdot \mathbf{n}) \boldsymbol{Q}^{\Xi_k}, \label{eq:Fcom_mortar}
\end{gather}
where the variables, without further explanation, have similar meanings to those in (\ref{eq:Fcom}).

\subsubsection{Project common values back to cell faces}
\label{sec:project_back}
\begin{figure}[!htb]
\centering
\includegraphics[scale=1.0]{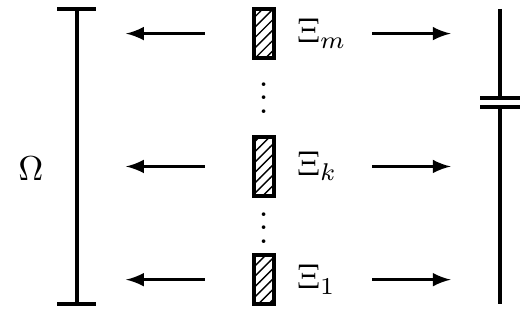}
\caption{Project common values from mortars back to cell faces.}
\label{fig:mortar_projection}
\end{figure}

Figure \ref{fig:mortar_projection} illustrates the process of projecting common values from $m$ mortars back to a cell face $\Omega$. Taking flux as an example, we can either directly project back the physical flux (method 1) or convert it to computational flux to project back (method 2). The details are described in what follows.

\pagebreak
\vspace{1.0em}
\noindent\textit{Method 1:}
\vspace{0.5em}

The inviscid normal fluxes on a cell face and a mortar are represented by the following polynomials,
\begin{equation}
F_\text{inv}^\Omega (\xi) = \sum_{i=1}^{N} F_{\text{inv},i}^\Omega h_i(\xi), \quad
F_\text{inv}^{\Xi_k} (z)  = \sum_{i=1}^{N} F_{\text{inv},i}^{\Xi_k} h_i(z).
\end{equation}
To get the $F_{\text{inv},i}^\Omega$'s, we require
\begin{equation}
\sum_{k=1}^{m} \int_{o_k}^{o_k+s_k} \left( F_\text{inv}^\Omega (\xi) - F_\text{inv}^{\Xi_k}(z) \right) h_j(\xi) ~\mathrm{d} \xi = 0, ~~ \forall j=1,2,...,N,
\label{eq:prjF0}
\end{equation}
which gives
\begin{equation}
\mathbf{M} \boldsymbol{F}_\text{inv}^\Omega = \sum_{k=1}^{m} s_k \mathbf{S}^{\Xi_k\rightarrow\Omega} \boldsymbol{F}_\text{inv}^{\Xi_k},
\label{eq:prjF1b0}
\end{equation}
where $\mathbf{M}$ is identical to that in (\ref{eq:prjQ0}), and $\mathbf{S}^{\Xi_k\rightarrow\Omega}$ is simply the transpose of the $\mathbf{S}^{\Omega\rightarrow\Xi_k}$ matrix from (\ref{eq:prjQ0}). Solutions of (\ref{eq:prjF1b0}) are
\begin{equation}
\boldsymbol{F}_\text{inv}^\Omega = \sum_{k=1}^{m} \mathbf{P}^{\Xi_k\rightarrow\Omega} \boldsymbol{F}_\text{inv}^{\Xi_k} = \sum_{k=1}^{m} s_k \mathbf{M}^{-1} \mathbf{S}^{\Xi_k\rightarrow\Omega} \boldsymbol{F}_\text{inv}^{\Xi_k},
\label{eq:prjF2b}
\end{equation}
where $\mathbf{P}^{\Xi_k\rightarrow\Omega}$ is the projection matrix from $\Xi_k$ to $\Omega$. The above physical flux on a cell face is then converted to a computational one to compute residuals. The conversion follows (\ref{eq:Fn2Ft}), here the normal is
\begin{equation}
\mathbf{N} = (y^\Omega_\xi,-x^\Omega_\xi) = L^\Omega (\cos\theta^\Omega, \sin\theta^\Omega),
\end{equation}
and the final computational flux on the cell face is
\begin{equation}
\widetilde{\boldsymbol{F}}\vphantom{\boldsymbol{F}}_\text{inv}^\Omega = \Vert \mathbf{N} \Vert \boldsymbol{F}_\text{inv}^\Omega = L^\Omega \boldsymbol{F}_\text{inv}^\Omega. \label{eq:Finv_back1}
\end{equation}

\vspace{1.0em}
\noindent\textit{Method 2:}
\vspace{0.5em}

Alternatively, we can convert the the physical flux in (\ref{eq:Fcom_mortar}) to computational, and then do the back projection. For a mortar, there are two types of normals depending on which space, i.e., the mortar space or the computational space, is taken as the reference space. For the mortar space, the normal is
\begin{equation}
\breve{\mathbf{N}} = (y^{\Xi_k}_z,-x^{\Xi_k}_z) = L^{\Xi_k} (\cos\theta^{\Xi_k}, \sin\theta^{\Xi_k}).
\end{equation}
For the computational space, the normal is
\begin{equation}
\widetilde{\mathbf{N}} = (y^{\Xi_k}_\xi,-x^{\Xi_k}_\xi) = (y^{\Xi_k}_z z_\xi,-x^{\Xi_k}_z z_\xi) = \frac{1}{s_k} (y^{\Xi_k}_z,-x^{\Xi_k}_z) = \frac{1}{s_k} \breve{\mathbf{N}}.
\end{equation}
The corresponding fluxes in these two spaces are
\begin{align}
\breve{\boldsymbol{F}}^{\Xi_k}_\text{inv} &= \Vert \breve{\mathbf{N}} \Vert \boldsymbol{F}_\text{inv}^{\Xi_k} = L^{\Xi_k} \boldsymbol{F}_\text{inv}^{\Xi_k}, \label{eq:Fxi1} \\[0.5em]
\widetilde{\boldsymbol{F}}\vphantom{\boldsymbol{F}}^{\Xi_k}_\text{inv} &= \Vert \widetilde{\mathbf{N}} \Vert \boldsymbol{F}_\text{inv}^{\Xi_k} = \frac{1}{s_k} \breve{\boldsymbol{F}}^{\Xi_k}_\text{inv}. \label{eq:Fcomp_mortar}
\end{align}
The computational inviscid fluxes on a cell face and a mortar are represented by the following polynomials,
\begin{equation}
\widetilde{F}_\text{inv}^{\Omega}(\xi) = \sum_{i=1}^{N} \widetilde{F}_{\text{inv},i}^{\Omega} h_i(\xi), \quad
\widetilde{F}_\text{inv}^{\Xi_k}(z)    = \sum_{i=1}^{N} \widetilde{F}_{\text{inv},i}^{\Xi_k}  h_i(z).
\end{equation}
To get the $\widetilde{F}_{\text{inv},i}^{\Omega}$'s, we require
\begin{equation}
\sum_{k=1}^{m} \int_{o_k}^{o_k+s_k} \left( \widetilde{F}_\text{inv}^{\Omega}(\xi) - \widetilde{F}_\text{inv}^{\Xi_k}(z) \right) h_j(\xi) ~\mathrm{d} \xi = 0, ~~ \forall j=1,2,...,N, \label{eq:prjF1}
\end{equation}
which leads to
\begin{equation}
\mathbf{M} \widetilde{\boldsymbol{F}}\vphantom{\boldsymbol{F}}_\text{inv}^{\Omega} = \sum_{k=1}^{m} s_k \mathbf{S}^{\Xi_k\rightarrow\Omega} \widetilde{\boldsymbol{F}}\vphantom{\boldsymbol{F}}_\text{inv}^{\Xi_k} = \sum_{k=1}^{m} \mathbf{S}^{\Xi_k\rightarrow\Omega} \breve{\boldsymbol{F}}_\text{inv}^{\Xi_k},
\label{eq:prjF1b}
\end{equation}
where $\mathbf{M}$ and $\mathbf{S}^{\Xi_k\rightarrow\Omega}$ are identical to those in (\ref{eq:prjF1b0}). The solutions of the above system are
\begin{equation}
\widetilde{\boldsymbol{F}}\vphantom{\boldsymbol{F}}_\text{inv}^{\Omega} = \sum_{k=1}^{m} \mathbf{P}^{\Xi_k\rightarrow\Omega} \breve{\boldsymbol{F}}_\text{inv}^{\Xi_k} = \sum_{k=1}^{m} \mathbf{M}^{-1} \mathbf{S}^{\Xi_k\rightarrow\Omega} \breve{\boldsymbol{F}}_\text{inv}^{\Xi_k},
\label{eq:prjF2}
\end{equation}
and note the difference between the $\mathbf{P}^{\Xi_k\rightarrow\Omega}$ matrices in the above equation and in Eq. (\ref{eq:prjF2b}).

These two methods are in fact equivalent as one may expect. To show this, simply divide both sides of (\ref{eq:prjF2}) by $L^\Omega$, and consider (\ref{eq:sk}), (\ref{eq:Finv_back1}), and (\ref{eq:Fxi1}), and we will get exactly the same result as (\ref{eq:prjF2b}). Also note the singularity in (\ref{eq:Fcomp_mortar}) when the scaling (size) of a mortar becomes zero. This singularity is eliminated in (\ref{eq:prjF1b}) where the scaling is multiplied back. Equation (\ref{eq:Fcomp_mortar}) is for derivation purpose only, and is not actually evaluated in the computation, so the singularity is naturally avoided.

\subsubsection{Treatment of viscous fluxes}
For viscous flow, we first project the common solution (\ref{eq:Qcom_mortar}) back to cell faces in the same way as (\ref{eq:prjF2b}). The updated solutions on a cell face are then involved in the computation of solution gradients. After that, we could either project local gradients (method 1) or local viscous fluxes (method 2) to mortars to compute common viscous fluxes which are then projected back. The steps are described below.

\vspace{1.0em}
\noindent\textit{Method 1:}
\vspace{0.5em}

The local gradients on cell faces are projected to mortars following (\ref{eq:prjQ2}). The common gradients and common physical viscous fluxes on a mortar are then computed following (\ref{eq:nablaq_com}) and (\ref{eq:Fvcom}), respectively. The normal viscous flux is calculated as
\begin{equation}
\breve{\boldsymbol{F}}^{\Xi_k}_\text{vis} = \accentset{\leftrightarrow}{\boldsymbol{F}}^{\Xi_k}_\text{vis} \cdot \breve{\mathbf{N}}, \label{eq:Fviscoms_mortar1}
\end{equation}
where $\accentset{\leftrightarrow}{\boldsymbol{F}}^{\Xi_k}_\text{vis} = (\boldsymbol{F}^{\Xi_k}_{\text{vis}}, \boldsymbol{G}^{\Xi_k}_{\text{vis}})$ with the two components representing the physical common viscous fluxes.  This normal viscous flux is finally projected back to cell faces following (\ref{eq:prjF2}).

\vspace{1.0em}
\noindent\textit{Method 2:}
\vspace{0.5em}

Local viscous flux, denoted by $\widetilde{\boldsymbol{F}}\vphantom{\boldsymbol{F}}_\text{vis}^\Omega$, is projected to mortars in the same way as (\ref{eq:prjQ1}). The resulting normal viscous flux on the left side of a mortar is
\begin{equation}
\breve{\boldsymbol{F}}_\text{vis}^{\Xi_k,\text{L}} = \mathbf{P}^{\Omega\rightarrow\Xi_k} \widetilde{\boldsymbol{F}}\vphantom{\boldsymbol{F}}_\text{vis}^\Omega = s_k \mathbf{M}^{-1} \mathbf{S}^{\Omega\rightarrow\Xi_k} \widetilde{\boldsymbol{F}}\vphantom{\boldsymbol{F}}_\text{vis}^\Omega. \label{eq:prjFvis}
\end{equation}
Note the difference on the $\mathbf{P}^{\Omega\rightarrow\Xi_k}$'s in the above equation and in (\ref{eq:prjQ2}). The viscous flux on the right side of a mortar is obtained in the same way. The common normal viscous flux is then calculated as
\begin{equation}
\breve{\boldsymbol{F}}^{\Xi_k}_\text{vis} =\frac{1}{2}(\breve{\boldsymbol{F}}^{\Xi_k,\text{L}}_\text{vis} + \breve{\boldsymbol{F}}^{\Xi_k,\text{R}}_\text{vis}), \label{eq:Fviscoms_mortar}
\end{equation}
which is projected back to cell faces following (\ref{eq:prjF2}) to update the original normal viscous flux.

We have compared method 1 and method 2 in a series of tests (not reported here), and do not notice any obvious difference on the results. But method 2 is slightly faster than method 1, because it requires fewer projections and calculations.

\subsection{Proof of Global Conservation}
\label{sec:conservation}
Global conservation means that, without the presence of source in a domain, the net gain or loss of $\widetilde{\mathbf{Q}}$ over the whole domain is determined solely by fluxes through the boundaries. Inspired by the works in \cite{kopriva-1996a,kopriva-1996b}, we prove in what follows that the present method is globally conservative.

Define the following quadrature weights $w_i$ and $w_j$ at the SPs such that
\begin{equation}
\int_{0}^{1} \int_{0}^{1} \phi(\xi,\eta) \text{d}\xi\text{d}\eta = \sum_{j=1}^{N} \sum_{i=1}^{N} \phi_{ij} w_i w_j, \quad \forall \phi \in \boldsymbol{\mathsf{P}}_{N-1,N-1}. \label{eq:wiwj}
\end{equation}
Multiply the discretized equation (\ref{eq:resid}) by the weights and sum over all the SPs
\begin{equation}
\sum_{j=1}^{N} \sum_{i=1}^{N} \left. \frac{\partial \widetilde{\mathbf{Q}}}{\partial t} \right|_{ij} w_i w_j = -\sum_{j=1}^{N} \sum_{i=1}^{N} \left[\frac{\partial \widehat{\mathbf{F}}}{\partial \xi} + \frac{\partial \widehat{\mathbf{G}}}{\partial \eta} \right]_{ij} w_i w_j. \label{eq:qfwiwj}
\end{equation}
Since the quadrature is exact for $\boldsymbol{\mathsf{P}}_{N-1,N-1}$ by its definition, and the terms $\partial \widetilde{\mathbf{Q}}/\partial t$, $\partial \widehat{\mathbf{F}}/\partial \xi$, $\partial \widehat{\mathbf{G}}/\partial \eta$ are in $\boldsymbol{\mathsf{P}}_{N-1,N-1}$, (\ref{eq:qfwiwj}) is thus equivalent to
\begin{equation}
\int_{0}^{1} \int_{0}^{1} \frac{\partial \widetilde{\mathbf{Q}}}{\partial t} \text{d}\xi\text{d}\eta = \int_{0}^{1} \int_{0}^{1} \left[\frac{\partial \widehat{\mathbf{F}}}{\partial \xi} + \frac{\partial \widehat{\mathbf{G}}}{\partial \eta} \right] \text{d}\xi\text{d}\eta.
\end{equation}
Take the time derivative out, and apply the divergence theorem to the right hand side,
\begin{equation}
\frac{\partial }{\partial t} \int_{0}^{1} \int_{0}^{1}  \widetilde{\mathbf{Q}} \text{d}\xi\text{d}\eta = \int_{0}^{1} \widehat{\mathbf{F}}(0,\eta) \text{d}\eta - \int_{0}^{1} \widehat{\mathbf{F}}(1,\eta) \text{d}\eta + \int_{0}^{1} \widehat{\mathbf{G}}(\xi,0) \text{d}\xi - \int_{0}^{1} \widehat{\mathbf{G}}(\xi,1) \text{d}\xi, \label{eq:cnsv_div}
\end{equation}
where, by definition, the corrected fluxes $\widehat{\mathbf{F}}$ and $\widehat{\mathbf{G}}$ take the common values on cell interfaces. For a conforming mesh, if we sum (\ref{eq:cnsv_div}) over all the mesh elements, the interior flux terms cancel on cell interfaces, leaving only the boundary flux terms. This confirms that the change of $\widetilde{\mathbf{Q}}$ over the whole domain depends only on boundary fluxes, i.e., the scheme is globally conservative on conforming mesh.

For nonconforming sliding mesh, we need to show that the total flux on the left side of a sliding interface equals that on the right side so that global conservation is satisfied. Since we have required (\ref{eq:prjF1}) to hold for both inviscid and viscous fluxes with any $h_j$, and the $h_j$'s form a basis of $\boldsymbol{\mathsf{P}}_{N-1}$, the following relation thus holds,
\begin{equation}
\sum_{k=1}^{m} \int_{o_k}^{o_k+s_k} \left( \widetilde{F}^{\Omega}(\xi) - \widetilde{F}^{\Xi_k}(z) \right) \psi(\xi) ~\mathrm{d} \xi = 0, \quad \forall \psi \in \boldsymbol{\mathsf{P}}_{N-1},
\end{equation}
Substitution of $\psi=1 \in \boldsymbol{\mathsf{P}}_{N-1}$ gives
\begin{equation}
\sum_{k=1}^{m} \int_{o_k}^{o_k+s_k} \left( \widetilde{F}^{\Omega}(\xi) - \widetilde{F}^{\Xi_k}(z) \right) \mathrm{d} \xi = 0.
\end{equation}
Rearrange the above equation to
\begin{equation}
\sum_{k=1}^{m} \int_{o_k}^{o_k+s_k} \widetilde{F}^{\Omega}(\xi) \mathrm{d}\xi = \sum_{k=1}^{m} \int_{o_k}^{o_k+s_k} \widetilde{F}^{\Xi_k}(z) \mathrm{d}\xi, \nonumber
\end{equation}
then change the variable of integration on the right hand side,
\begin{equation}
\int_{0}^{1} \widetilde{F}^{\Omega}(\xi) \mathrm{d} \xi = \sum_{k=1}^{m} \int_{0}^{1} s_k \widetilde{F}^{\Xi_k}(z) \mathrm{d} z, \nonumber
\end{equation}
and now consider the relation in (\ref{eq:Fcomp_mortar}) (which also applies to viscous flux),
\begin{align}
\int_{0}^{1} \widetilde{F}^{\Omega}(\xi) \mathrm{d} \xi = \sum_{k=1}^{m} \int_{0}^{1} \breve{F}^{\Xi_k}(z) \mathrm{d} z. \label{eq:cnsv_omg}
\end{align}
The relation in (\ref{eq:cnsv_omg}) is valid for any sliding cell face $\Omega$, and it says that the flux on a sliding cell face equals the total flux on all its mortars. Since each mortar is unique and does not overlap with other mortars (see Fig. \ref{fig:mortar_schematic}), summing (\ref{eq:cnsv_omg}) over all the cell faces on the left side of a sliding interface will lead us to the conclusion that the total flux on the left side of a sliding interface equals the total flux on all the mortars. The same conclusion also holds on the right side of a sliding interface. Therefore, the total flux on the left side of a sliding interface equals that on the right side, i.e., the scheme satisfies global conservation.

\subsection{Proof of Outflow Condition}
In a hyperbolic system, the propagation of waves should not be affected by downwind signals, which is called the outflow condition. When it comes to our case, satisfying the outflow condition means that if we project a variable from a cell face to its mortars, and then immediately project back, the value of the variable should not change. In what follows, we prove that the present sliding-mesh method satisfies this condition. Since the two projection methods in Sec. \ref{sec:project_back} are equivalent, we only show the poof on the first method for simplicity.

Let $\phi^\Omega$ denote a variable on a cell face, $\phi^{\Xi_k}$ the projected variable on the $k$-th mortar, and $\phi^{*\Omega}$ the variable from back projection. Then, the problem becomes proving that $\phi^{*\Omega}=\phi^\Omega$.

Start with the following polynomial representations of these variables,
\begin{equation}
\phi^{\Omega} (\xi) = \sum_{i=1}^{N} \phi_i^{\Omega}h_i(\xi), \quad \phi^{\Xi_k} (z) = \sum_{i=1}^{N} \phi_i^{\Xi_k}h_i(z), \quad \phi^{*\Omega} (\xi) = \sum_{i=1}^{N} \phi_i^{*\Omega}h_i(\xi).
\end{equation}
The projection from cell face to mortar follows (\ref{eq:prjQ1}),
\begin{equation}
\int_0^1 \left( \phi^{\Xi_k}(z) - \phi^{\Omega}(\xi) \right) h_j(z) ~\mathrm{d}z = 0, \quad \forall j=1,2,...,N,  \label{eq:phiz}
\end{equation}
from which we can get the projected variable on the mortar,
\begin{equation}
\boldsymbol{\phi}^{\Xi_k} = \mathbf{P}^{\Omega\rightarrow\Xi_k} \boldsymbol{\phi}^{\Omega} = \mathbf{M}^{-1} \mathbf{S}^{\Omega\rightarrow\Xi_k} \boldsymbol{\phi}^{\Omega}. \label{eq:phivec}
\end{equation}
Now, project $\phi^{\Xi_k}$ back to the cell face following (\ref{eq:prjF0}),
\begin{equation}
\sum_{k=1}^{m} \int_{o_k}^{o_k+s_k} \left( \phi^{*\Omega} (\xi) - \phi^{\Xi_k}(z) \right) h_j(\xi) ~\mathrm{d} \xi = 0, \quad \forall j=1,2,...,N, \label{eq:phistar}
\end{equation}
from which we can get the variable from back projection,
\begin{equation}
\boldsymbol{\phi}^{*\Omega} = \sum_{k=1}^{m} \mathbf{P}^{\Xi_k\rightarrow\Omega} \boldsymbol{\phi}^{\Xi_k} = \sum_{k=1}^{m} s_k \mathbf{M}^{-1} \mathbf{S}^{\Xi_k\rightarrow\Omega} \boldsymbol{\phi}^{\Xi_k}. \label{eq:phisvec}
\end{equation}
Since the the $h_j$'s form a basis for $\boldsymbol{\mathsf{P}}_{N-1}$, from (\ref{eq:phiz}) we can get
\begin{equation}
\int_0^1 \left( \phi^{\Xi_k}(z) - \phi^{\Omega}(\xi) \right) \psi ~\mathrm{d}z = 0, \quad \forall \psi \in \boldsymbol{\mathsf{P}}_{N-1}.
\end{equation}
Because $h_j(\xi)=h_j(o_k+s_k z) \in \boldsymbol{\mathsf{P}}_{N-1}$, setting $\psi=h_j(\xi)$ gives
\begin{equation}
\int_0^1 \left( \phi^{\Xi_k}(z) - \phi^{\Omega}(\xi) \right) h_j(\xi) ~\mathrm{d}z = 0, \nonumber
\end{equation}
now rearrange the terms and change the variable of integration,
\begin{align}
\int_{o_k}^{o_k+s_k} \phi^{\Omega}(\xi) h_j(\xi) ~\mathrm{d}\xi = s_k \int_{0}^{1} \phi^{\Xi_k}(z) h_j(\xi) ~\mathrm{d}z. \label{eq:ofc1}
\end{align}
Apply (\ref{eq:ofc1}) to all the mortars of $\Omega$ and sum them up,
\begin{gather}
\sum_{k=1}^{m} \int_{o_k}^{o_k+s_k} \phi^{\Omega}(\xi) h_j(\xi) ~\mathrm{d}\xi = \sum_{k=1}^{m} s_k \int_{0}^{1} \phi^{\Xi_k}(z) h_j(\xi) ~\mathrm{d}z, \nonumber \\[0.5em]
\int_{0}^{1} \phi^{\Omega}(\xi) h_j(\xi) ~\mathrm{d}\xi = \sum_{k=1}^{m} s_k \int_{0}^{1} \phi^{\Xi_k}(z) h_j(\xi) ~\mathrm{d}z, \label{eq:phixixi1}
\end{gather}
On the other hand,  rearranging (\ref{eq:phistar}) and changing the variable of integration gives
\begin{align}
\int_{0}^{1} \phi^{*\Omega} (\xi) h_j(\xi) \mathrm{d}\xi = \sum_{k=1}^{m} s_k \int_{0}^{1} \phi^{\Xi_k}(z) h_j(\xi) ~\mathrm{d}z. \label{eq:phistar2}
\end{align}
Note that (\ref{eq:phixixi1}) and (\ref{eq:phistar2}) have exactly the same right hand side, therefore
\begin{align}
\int_{0}^{1} \phi^{*\Omega} (\xi) h_j(\xi) \mathrm{d}\xi = \int_{0}^{1} \phi^{\Omega} (\xi) h_j(\xi) \mathrm{d}\xi, \quad \forall j=1,2,...,N,
\end{align}
which is equivalent to the following equation system,
\begin{equation}
\mathbf{M} \boldsymbol{\phi}^{*\Omega} = \mathbf{M} \boldsymbol{\phi}^{\Omega}. \nonumber
\end{equation}
Because $\mathbf{M}$ is invertible, we therefore have
\begin{equation}
\boldsymbol{\phi}^{*\Omega} = \mathbf{M}^{-1}\mathbf{M} \boldsymbol{\phi}^{\Omega}, \nonumber
\end{equation}
that is
\begin{equation}
\boldsymbol{\phi}^{*\Omega} = \boldsymbol{\phi}^{\Omega}, \label{eq:ofc2}
\end{equation}
which is exactly what we were expecting, and the present method therefore satisfies the outflow condition. If we substitute (\ref{eq:phivec}) into (\ref{eq:phisvec}) and take (\ref{eq:ofc2}) into account, we have the following relation
\begin{equation}
\boldsymbol{\phi}^{*\Omega} = \sum_{k=1}^{m} \mathbf{P}^{\Xi_k\rightarrow\Omega} \mathbf{P}^{\Omega\rightarrow\Xi_k} \boldsymbol{\phi}^{\Omega} = \boldsymbol{\phi}^{\Omega},
\end{equation}
which readily reveals a property of the projection matrices,
\begin{equation}
\sum_{k=1}^{m} \mathbf{P}^{\Xi_k\rightarrow\Omega} \mathbf{P}^{\Omega\rightarrow\Xi_k} = \mathbf{I},
\label{eq:outflow}
\end{equation}
where $\mathbf{I}$ is the identity matrix.

We have numerically tested (\ref{eq:outflow}) on many cases, and it holds perfectly without exception. As it was noted in \cite{kopriva-1996b} that satisfaction of the outflow condition ensures that a method does not change the characteristics of the flow field, and therefore should not degrade the maximum allowable time step size compared with the original method on conforming mesh. Our observation through numerical tests is consistent with this expectation.

\subsection{On the Implementation}
\label{subsec:implement}

\begin{figure}[!htb]
\centering
\includegraphics[scale=1.0]{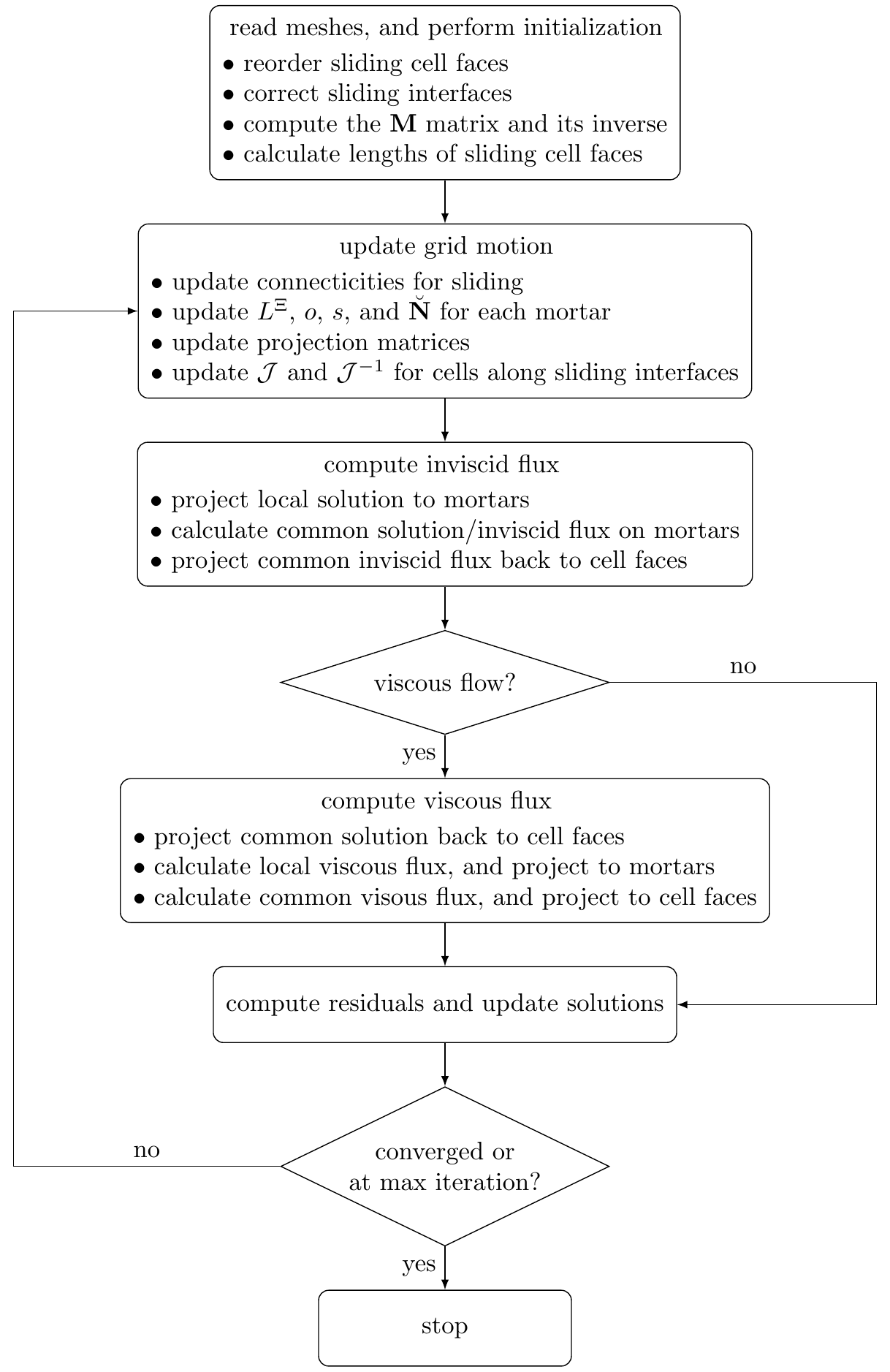}
\caption{Flow chart of the implementation procedures.}
\label{fig:flowchart}
\end{figure}

A flow chart of the implementation procedures is shown in Fig. \ref{fig:flowchart}. Some of the steps are explained with more details in what follows.

\subsubsection{Read meshes}
\begin{figure}[!htb]
\centering
\includegraphics[scale=1.0]{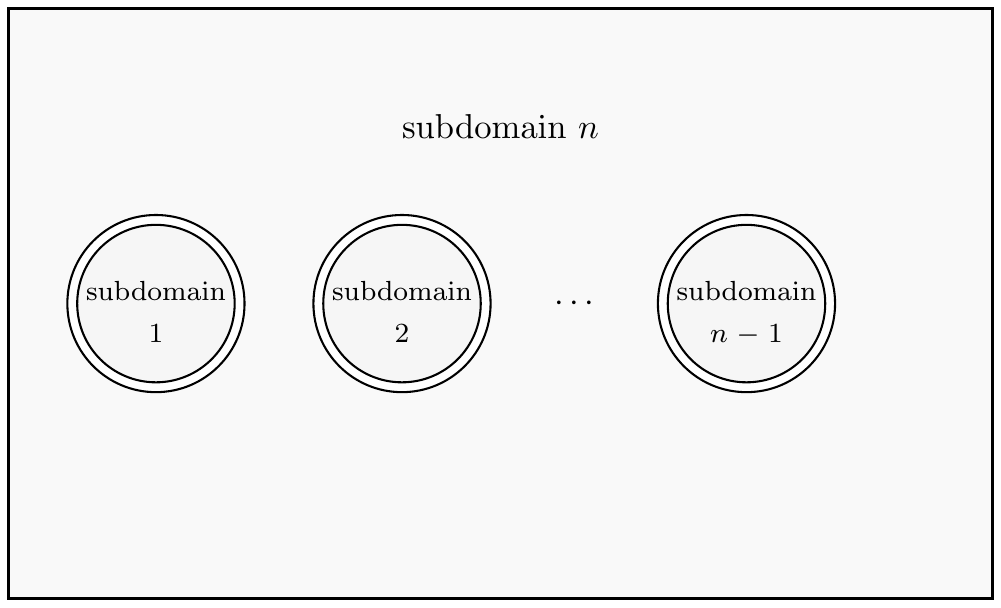}
\caption{Schematic of $n$ subdomains with $(n-1)$ sliding interfaces (inner subdomains have been scaled).}
\label{fig:subdomains}
\end{figure}

Mesh for each subdomain is generated independently and stored in a separate file. When subdomain meshes are read in, they are assembled into a ``single'' mesh by adding offsets to the numberings of cells and vertices. For example, as shown in Fig. \ref{fig:subdomains}, assume we have $n$ subdomains with $N_1$, $N_2$, ..., $N_n$ cells. In the assembled mesh, numbering for the cells in subdomain 1 starts from $1$, for subdomain 2 starts from $N_1+1$, for subdomain 3 starts from $N_1+N_2+1$, and so on. The same rule applies to vertices and boundary faces. In this way, cells and vertices are numbered uniquely and continuously. The overall mesh behaves just like a single mesh and sliding interfaces are like interior boundaries.

\subsubsection{Reorder sliding cell faces}
All sliding cell faces are stored consecutively in an array. This means that cell faces from the left side of the first sliding interface are stored first, followed by those from the right side of the first sliding interface, and then those from the left side of the second sliding interface, and so on. Cell faces from each side are reordered into counterclockwise order, which makes connectivity updating more efficient. The algorithm for the reordering is quite simple: any cell face can be assigned as the first cell face, and then the next cell face is the one whose starting vertex is the ending vortex of the current cell face, repeat until finished.

\subsubsection{Correct sliding interface}
Mesh generator may introduce geometric errors. When it comes to sliding mesh, the problem is usually that, the vertices on a circular sliding interface do not represent the same radius. This geometric error may potentially contaminate a simulation. Here is an easily fix that can be performed during preprocessing: pick up a vertex and take the corresponding radius as a reference radius, and then use this reference radius and the original angle of each vertex to update the coordinates of that vertex. In this way, all the vertices on a circular sliding interface represent the same radius.

\subsubsection{Update connectivities}
For simplicity, we only take the first sliding interface from Fig. \ref{fig:subdomains} to explain how the connectivities are efficiently updated. The same procedure is repeated for other sliding interfaces (there is no limit on how many we may have). Table \ref{tab:mortar_connectivity} shows the definitions of important numbers and arrays used during this process. The detailed procedures for updating the connectivities are listed as Algorithm \ref{alg:mortar}. The general idea is that we ``walk'' along each sliding interface counterclockwise, and two successive points form a mortar (these two points might or might not come from the same side). One may also refer to Fig. \ref{fig:mortar_schematic} to facilitate the understanding of the algorithm.
\begin{table}[!htb]
\def\arraystretch{1.4}
\setlength{\tabcolsep}{5mm}
\centering
\begin{tabular}{|>{\centering}m{0.14\textwidth} | m{0.58\textwidth}|}
\hline
\multicolumn{1}{|c|}{Variable}  &  \multicolumn{1}{c|}{Meaning} \\
\hline
\textit{nfl} & number of cell faces on the left side   \\
\hline
\textit{nfr}  & number of cell faces on the right side  \\
\hline
\textit{nf}  & total number of faces on this sliding interface  \newline note: $\mathit{nf=nfl+nfr}$ (by definition) \\
\hline
\textit{nm}  & total number of mortars on this sliding interface \newline note: $\mathit{nm=nf}$ (always valid)  \\
\hline
\texttt{vof}(1:2,1:\textit{nf})  & \texttt{vof}(1:2,\textit{ifa}) stores the two vertices of the \textit{ifa}-th cell face \\
\hline
\texttt{mof}(1:2,1:\textit{nf})  & $\texttt{mof}(1,\mathit{ifa})$ stores the first mortar of the \textit{ifa}-th cell face \newline $\texttt{mof}(2,\mathit{ifa})$ stores the number of mortars of the \textit{ifa}-th cell face  \\
\hline
\texttt{fom}(1:2,1:\textit{nm})  & $\texttt{fom}(1,\mathit{im})$ stores the left cell face of the \textit{im}-th mortar \newline $\texttt{fom}(2,\mathit{im})$ stores the right cell face of the \textit{im}-th mortar  \\
\hline
\texttt{vom}(1:2,1:\textit{nm})  & \texttt{vom}(1:2,\textit{im}) stores the two vertices of the \textit{im}-th mortar  \\
\hline
\end{tabular}
\caption{Definitions of variables used in updating cell face and mortar connectivities on a sliding interface.}
\label{tab:mortar_connectivity}
\end{table}

\begin{algorithm}[!htb]
\caption{Algorithm for updating cell face and mortar connectivities for a sliding interface} \label{alg:mortar}
\algnewcommand\algorithmicto{\textbf{to}}
\algrenewtext{For}[2]{\algorithmicfor\ #1 \algorithmicto\ #2 \algorithmicdo}
\begin{algorithmic}[0]
\State \texttt{mof}=0; \texttt{fom}=0; \texttt{vom}=0  \Comment initialize with zeros
\State
\State $\mathit{ifl}=1$
\For{$\mathit{ifr=(nfl}+1)$}{$\mathit{nf}$}
\If{$\texttt{vof}(1,\mathit{ifl})$ lies between $\texttt{vof}(1,\mathit{ifr})$ and $\texttt{vof}(2,\mathit{ifr})$}  \Comment the first mortar is located
\State exit
\EndIf
\EndFor
\State
\State $\mathit{im} \gets 1$
\State $\texttt{mof}(1,\mathit{ifl}) \gets \mathit{im}$                        \Comment connectivities of the first mortar
\State $\texttt{mof}(2,\mathit{ifl}) \gets \texttt{mof}(2,\mathit{ifl}) + 1$
\State $\texttt{mof}(2,\mathit{ifr}) \gets \texttt{mof}(2,\mathit{ifr}) + 1$
\State $\texttt{fom}(1,\mathit{im}) \gets \mathit{ifl}$
\State $\texttt{fom}(2,\mathit{im}) \gets \mathit{ifr}$
\State $\texttt{vom}(1,\mathit{im}) \gets \texttt{vof}(1,\mathit{ifl})$
\State
\For{$\mathit{im}=2$}{$\mathit{nm}$}      \Comment connectivities of remaining mortars
\If{$\texttt{vof}(2,\mathit{ifl})$ lies between $\texttt{vof}(1,\mathit{ifr})$ and $\texttt{vof}(2,\mathit{ifr})$}
\State $\mathit{ifl} \gets \mathit{ifl} + 1$
\State $\mathit{ifa} \gets \mathit{ifl}$
\Else
\State $\mathit{ifr} \gets \mathit{ifr} + 1$
\If{$(\mathit{ifr>nf})$}
\State $\mathit{ifr} \gets \mathit{ifr-nfr}$
\EndIf
\State $\mathit{ifa} \gets \mathit{ifr}$
\EndIf
\State
\State $\texttt{mof}(1,\mathit{ifa}) \gets \mathit{im}$
\State $\texttt{mof}(2,\mathit{ifl}) \gets \texttt{mof}(2,\mathit{ifl}) + 1$
\State $\texttt{mof}(2,\mathit{ifr}) \gets \texttt{mof}(2,\mathit{ifr}) + 1$
\State $\texttt{fom}(1,\mathit{im}) \gets \mathit{ifl}$
\State $\texttt{fom}(2,\mathit{im}) \gets \mathit{ifr}$
\State $\texttt{vom}(1,\mathit{im}) \gets \texttt{vof}(\mathit{1,ifa})$
\State $\texttt{vom}(2,\mathit{im}-1) \gets \texttt{vof}(\mathit{1,ifa})$
\EndFor
\State
\State $\texttt{vom}(2,\mathit{nm}) \gets \texttt{vom}(1,1)$   \Comment mortars always form a closed loop
\end{algorithmic}
\end{algorithm}

\subsubsection{Compute fluxes}
The flow chart in Fig. \ref{fig:flowchart} only shows the steps related to sliding. These steps must be injected into the main solver at the right locations. For example, the ``compute inviscid flux'' part is called after interior and cell interface inviscid flux computations are done. The back projection of common solution takes place after all inviscid flux computations are finished and before viscous flux computations start. The projections of viscous flux are carried out after all interior and cell interface viscous flux computations are done.

\section{Examples}
\label{sec:examples}
In this section, we apply the sliding mesh method to several flow problems. We first test it on an inviscid flow and a viscous flow to verify the spatial and temporal accuracies. We then study the conservation and the free-stream-preservation properties of the method. Following that, we perform a comparison study between the present method and a rigid-rotation method on flow over a rotating square cylinder. Finally, we apply the method to simulate flow over multiple rotating square cylinders to further demonstrate its capability.

For all the test cases wherever applicable, we employ the following  $L_2$ norm to measure the errors,
\begin{equation}
L_2 \text{ error} = \sqrt{ \frac{\sum_{i=1}^{N_\text{DOF}} (\phi_i - \phi_i^{\mathrm{exact}})^2}{N_\text{DOF}}},
\end{equation}
where $\phi$ represents the variable of interest, $\phi_i$ and $\phi_i^\text{exact}$ are the numerical and the exact solutions at the $i$-th degree of freedom (DOF), $N_\text{DOF}=N_\text{elem} \cdot N^2$ is the total number of DOFs, $N_\text{elem}$ is the total number of mesh elements, $N=P+1$ is the scheme order (also the number of SPs in each direction of an element), and $P$ is the polynomial degree.

\subsection{Spatial Accuracy}

\subsubsection{Euler vortex flow}
This is an inviscid flow test. In this flow, an isentropic vortex is superimposed to and convected by a uniform flow. The flow field in an infinite domain at a time instant $t$ can be analytically expressed as
\begin{align}
u    &= U_{\infty}   \left\{\cos\theta - \frac{\epsilon y_r}{r_c}\exp\left(\frac{1-x_r^2-y_r^2}{2r_c^2}\right) \right\}, \\[0.0em]
v    &= U_{\infty}   \left\{\sin\theta + \frac{\epsilon x_r}{r_c}\exp\left(\frac{1-x_r^2-y_r^2}{2r_c^2}\right) \right\}, \\[0.0em]
\rho &= \rho_{\infty}\left\{1 - \frac{(\gamma-1) (\epsilon M_{\infty})^2}{2} \exp\left(\frac{1-x_r^2-y_r^2}{r_c^2}\right) \right\}^{\frac{1}{\gamma -1}}, \\[0.0em]
p    &= p_{\infty}   \left\{1 - \frac{(\gamma-1) (\epsilon M_{\infty})^2}{2} \exp\left(\frac{1-x_r^2-y_r^2}{r_c^2}\right) \right\}^{\frac{\gamma}{\gamma -1}},
\end{align}
where $U_{\infty}$, $\rho_{\infty}$, $p_{\infty}$, $M_{\infty}$ are the speed, density, pressure and Mach number of the uniform flow; $\theta$ denotes the mean flow direction; $\epsilon$ and $r_c$ are the strength and size of the vortex; $(x_r,y_r)=(x - x_0 - \bar{u} t,y - y_0 - \bar{v} t)$ are the relative coordinates; $(x_0,y_0)$ represent the initial position of the vortex;  $(\bar{u},\bar{v})=(U_{\infty}\cos\theta,U_{\infty}\sin\theta)$ are the velocity components of the mean flow.

\begin{figure}[!htb]
\centering
\includegraphics[scale=0.44]{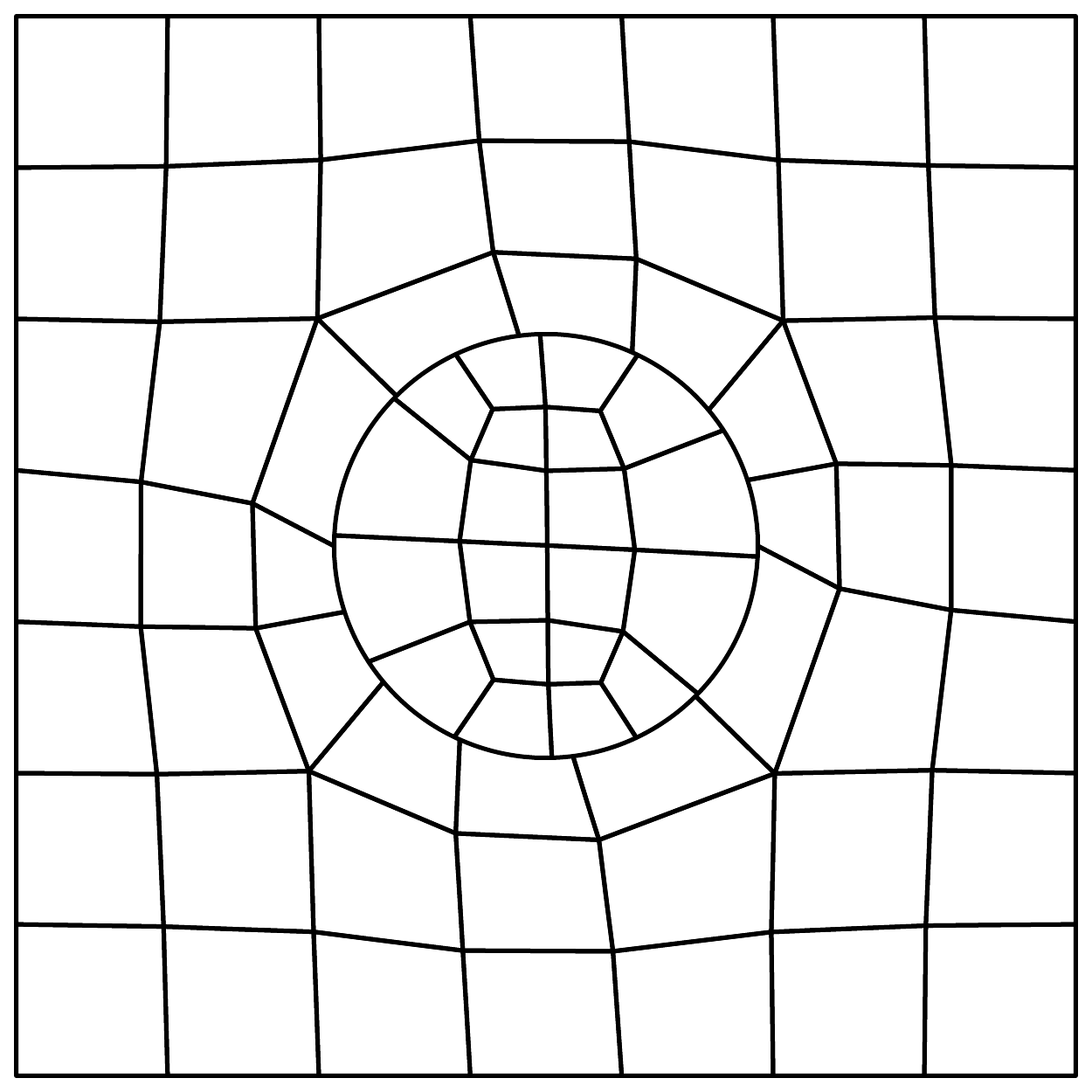}
\caption{Mesh for Euler vortex flow simulation.}
\label{fig:euler_msh}
\end{figure}

For this test, we have chosen the following parameters: $U_\infty = 1$, $\rho_\infty = 1$, $M_{\infty}=0.3$, $\theta=\arctan(1/2)$, $\epsilon=1$, and $r_c=1$. The overall computational domain has a size of $0\leq x, y \leq 10$, and the vortex is initially placed at $(x_0,y_0) = (5,5)$. Time-dependent analytical solutions are weakly prescribed along the boundaries to provide Dirichlet boundary conditions. The mesh used for this simulation is shown in Fig. \ref{fig:euler_msh}, where there are 72 mesh elements in total, with 20 of them in the inner rotating subdomain whose radius is $2$. Six rotational speeds: $\omega=0,1,5,10$, $15$ and $20$, are tested to study their effects on the solution. For all the cases, a five-stage fourth-order SSP Runge-Kutta scheme \cite{kraaijevanger-1991,ruuth-2006} with a time step size of $\Delta t=1.0\times 10^{-4}$ is used for the time marching.

\begin{figure}[!htb]
\centering
\includegraphics[width=\textwidth]{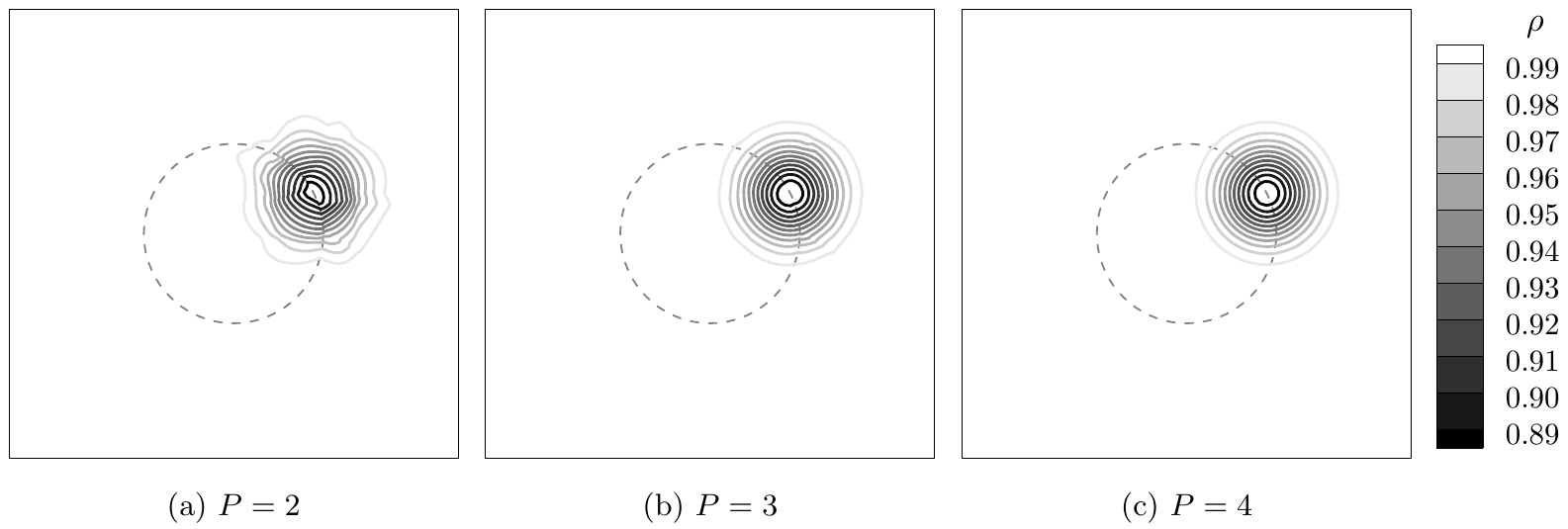}
\caption{Density contours of Euler vortex flow at $t=2$ (dashed lines represent sliding interface).}
\label{fig:euler_rho}
\end{figure}

In Fig. \ref{fig:euler_rho}, we compare the density contours at $t=2$ from the $\omega=20$ case using different polynomials. At this time instant, the vortex center travels right onto the sliding interface. It is obvious that $P=2$ does not provide enough resolution as the vortex is poorly resolved. But as the polynomial degree increases, the solution quality becomes much improved. Even at a small polynomial degree of $P=4$ and on such a coarse grid, the details of the vortex are very well captured.

\begin{figure}[!htb]
\centering
\includegraphics[scale=1.0]{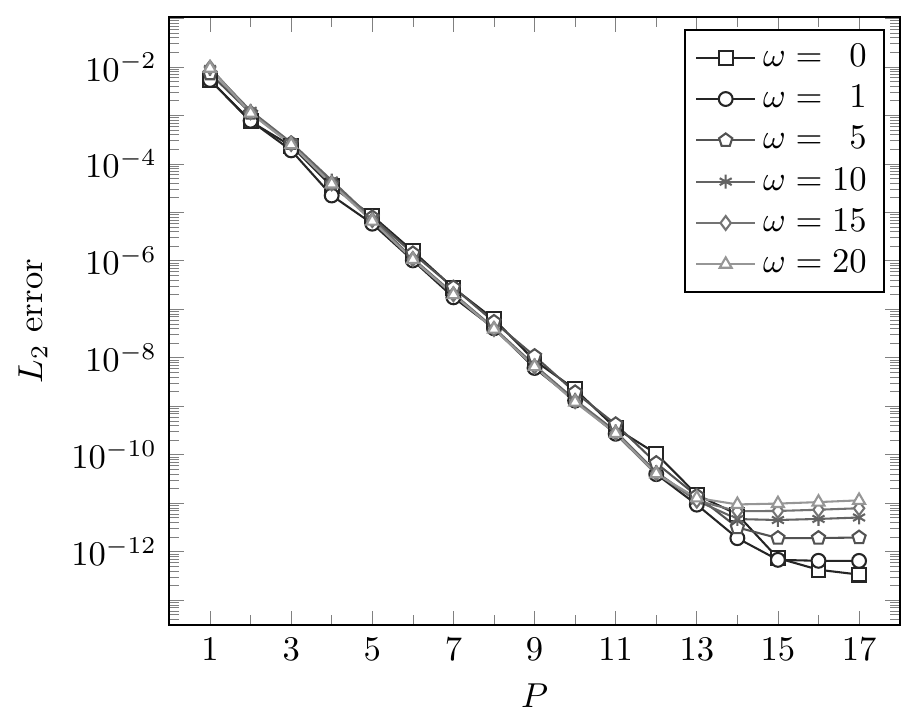}
\caption{$L_2$ error of density against polynomial degree for Euler-vortex flow simulation.}
\label{fig:euler_order}
\end{figure}

The $L_2$ errors of density are plotted in Fig. \ref{fig:euler_order} against polynomial degrees. It is seen that the errors from different cases are comparable when $P\le 13$, and start differing when $P>13$. The reason for this is that spatial errors dominate in the first regime, and they keep decreasing to such a small level that temporal errors start to dominate in the second regime. And larger rotational speed induces larger temporal errors as the second regime of Fig. \ref{fig:euler_order} shows. Similar observation on the effects of rotational speed was reported in \cite{ferrer-2012}. But in that work, the effects start showing up at a very early stage of $P=4$, which is mainly due to the low-order (2nd order) temporal scheme used in that work. The present results suggest that the rotational speed effects can be dramatically reduced by using a high-order temporal scheme (4th order here).

A closer look at the curves in Fig. \ref{fig:euler_order} reveals ``inconsistencies'' on the results. For example, the $\omega=20$ case generally gives smaller errors than the $\omega=5$ case does in the first regime, and the errors are sometimes even smaller than those of the $\omega=0$ case. These inconsistencies originate from the nonuniformity of the mesh in the rotating subdomain. When the vortex travels to the same location, it is actually captured by different grid resolutions when the rotational speed differs. Nevertheless, for all the cases, the errors inevitably decrease exponentially, which confirms the high-order accuracy of the sliding-mesh method.

\subsubsection{Taylor-Couette flow}
Taylor-Couette flow is formed between two concentric rotating circular cylinders. Due to viscous effects, this flow will reach a steady state if the Reynolds number is small. The steady-state azimuthal flow speed has the following expression,
\begin{equation}
v_{\theta} = v_{\theta_i} \frac{r_o/r - r/r_o}{r_o/r_i - r_i/r_o} + v_{\theta_o} \frac{r/r_i - r_i/r}{r_o/r_i - r_i/r_o},
\label{eq:couete_v}
\end{equation}
where $r_i$ and $r_o$ are the radii of the inner and the outer boundaries; $v_{\theta_i}$ and $v_{\theta_o}$ are the azimuthal flow speed at these two boundaries.

\begin{figure}[!htb]
\centering
\includegraphics[scale=0.5]{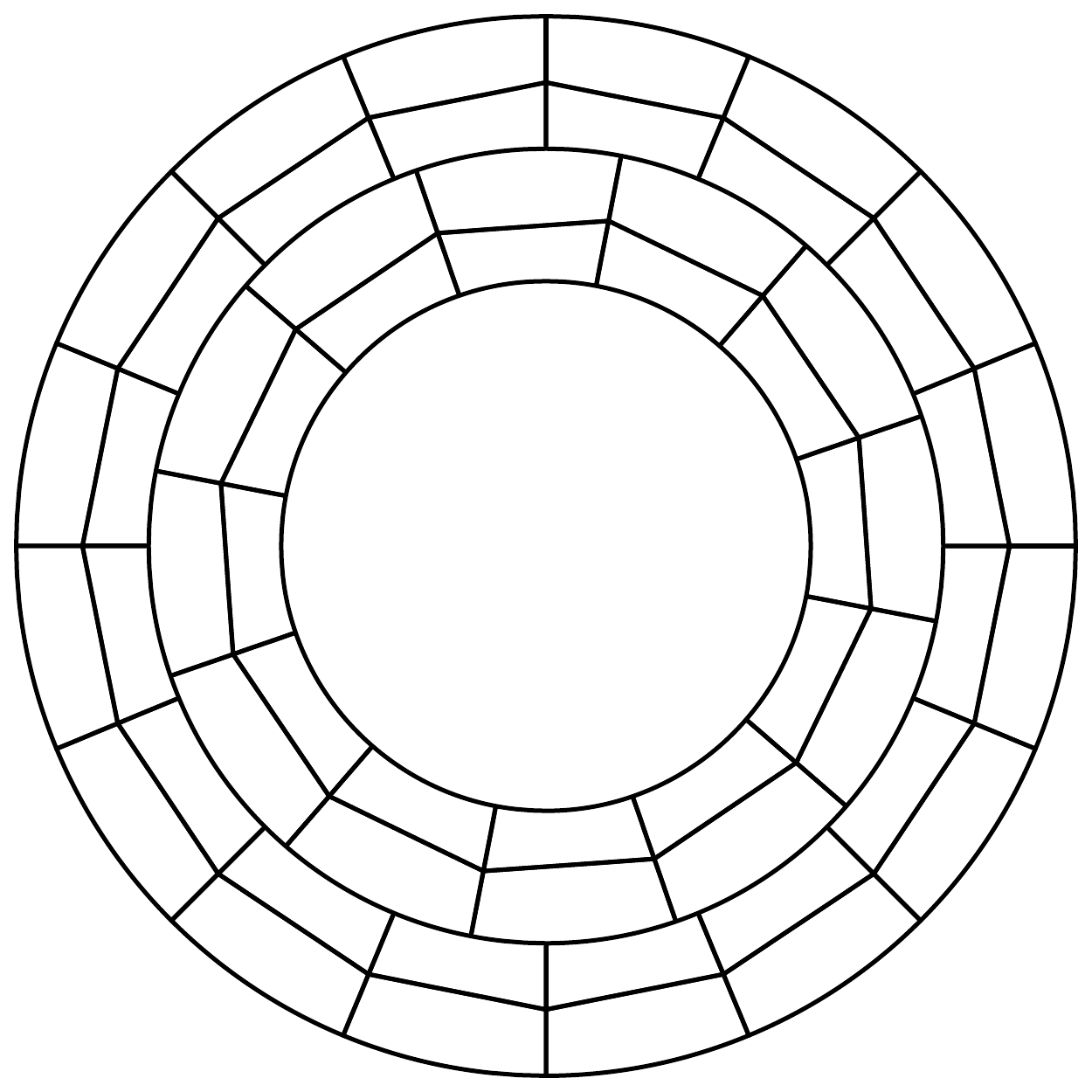}
\caption{Mesh for Taylor-Couette flow simulation.}
\label{fig:couette_msh}
\end{figure}

Figure \ref{fig:couette_msh} shows the mesh for this simulation. The overall domain is bounded by $r_i=1$ and $r_o=2$, and is split into a rotating inner subdomain and a fixed outer subdomain by a sliding interface at $r_s=1.5$. These two subdomains are meshed into 24 and 32 elements, with uniform grid distribution in azimuthal and radial directions.
The outer boundary is treated as a no-slip isothermal wall with $v_{\theta_o}=0$. The inner boundary is set as a Dirichlet boundary with $\rho_i=1$, Mach number $\mathit{Ma}=0.1$, and $v_{\theta_i}=1$. The Reynolds number based on flow properties at the inner boundary is $Re=10$. Again, six rotational speeds: $\omega=0,1,5,10$, $15$ and $20$, are tested. The same time marching scheme as that for the previous test with a time step size of $\Delta t=5.0\times 10^{-5}$ is used for all the tests. Exact transfinite mapping is employed on all circular boundaries, and all interior cell faces are represented linearly.

\begin{figure}[!htb]
\centering
\includegraphics[width=\textwidth]{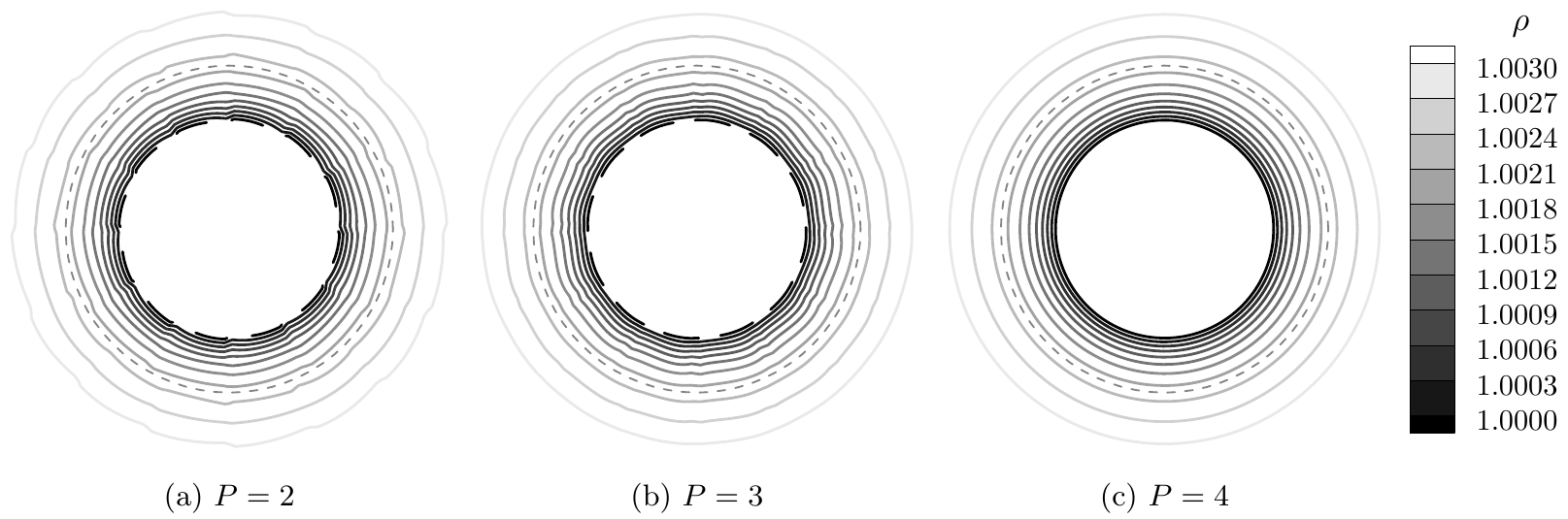}
\caption{Steady state density contours of Taylor-Couette flow (dashed lines in the middle represent sliding interface).}
\label{fig:couette_contour}
\end{figure}

Figure \ref{fig:couette_contour} shows the steady state density contours from different polynomials. The contours are expected to be a series of concentric circles, and we see improving results as the polynomial degree increases. Note that since the boundary conditions are weakly imposed, the computed values on the inner boundary are therefore not necessarily identical to the prescribed values. And this is the reason for the dashed lines at the inner boundaries. But the boundary conditions becomes stricter as the polynomial degree increases, which is why the dashed lines have disappeared at $P=4$.

\begin{figure}[!hbt]
\centering
\includegraphics[scale=1.0]{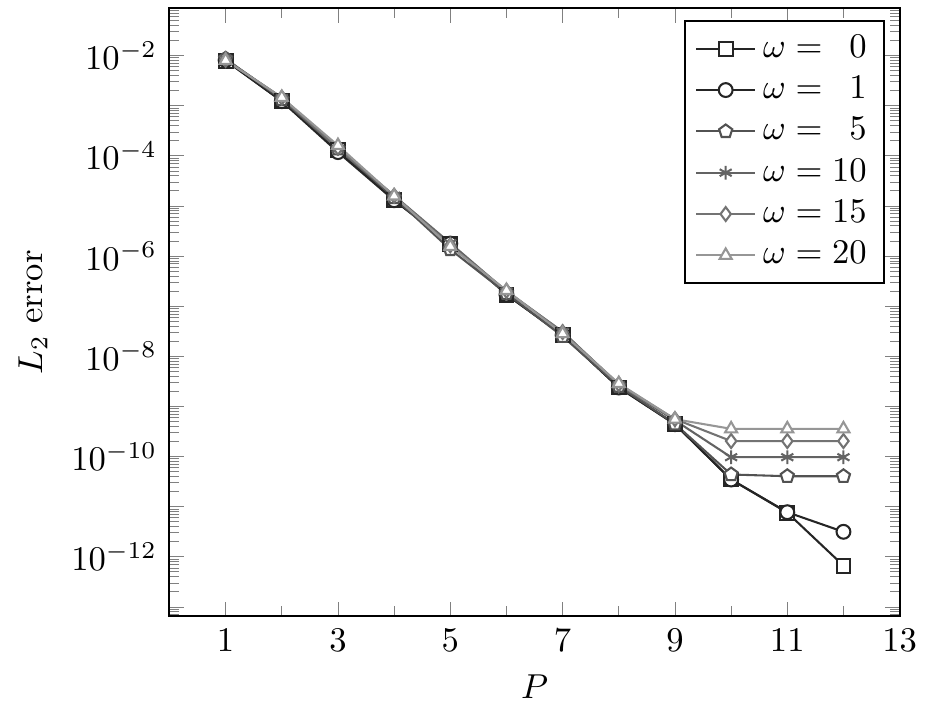}
\caption{$L_2$ errors of $u$ velocity against polynomial degree for Taylor-Couette flow.}
\label{fig:couette_order}
\end{figure}

We compute the $L_2$ errors of the $u$ velocity (i.e., the $x$-component of $v_\theta$) and plot the results in Fig. \ref{fig:couette_order}. Overall, the errors decrease exponentially with polynomial degree before temporal errors become dominant, which confirms the high-order accuracy of the present method on viscous flow. Meanwhile, the results from different rotational speeds are more consistent than those from the previous test. This is because of the uniformity of the current mesh in the azimuthal direction, which gives identical mesh resolution even when the rotational speed differs. The effects of rotational speed again  agree with our expectation.

\subsection{Temporal Accuracy}
In this section, we test the temporal accuracy and verify that the present sliding-mesh method does not affect the order of accuracy of a time marching scheme. We employ several strong-stability-preserving Runge-Kutta (SSPRK) schemes for the temporal discretization. Following the rules in \cite{spiteri-2002}, we denote an $s$-stage $p$-th order SSPRK scheme as SSP($s,p$). The CFL condition requires small enough time step size to stabilize a scheme, which makes it difficult for temporal error to dominate. To alleviate this issue, we employ the following SSPRK schemes that allow larger CFL numbers (and thus larger time step sizes): the SSP(4,2) scheme from \cite{spiteri-2002}, the SSP(8,3) scheme from \cite{ruuth-2006}, and the SSP(10,4) scheme from \cite{ketcheson-2008}.

\begin{figure}[!hbt]
\centering
\includegraphics[width=\textwidth]{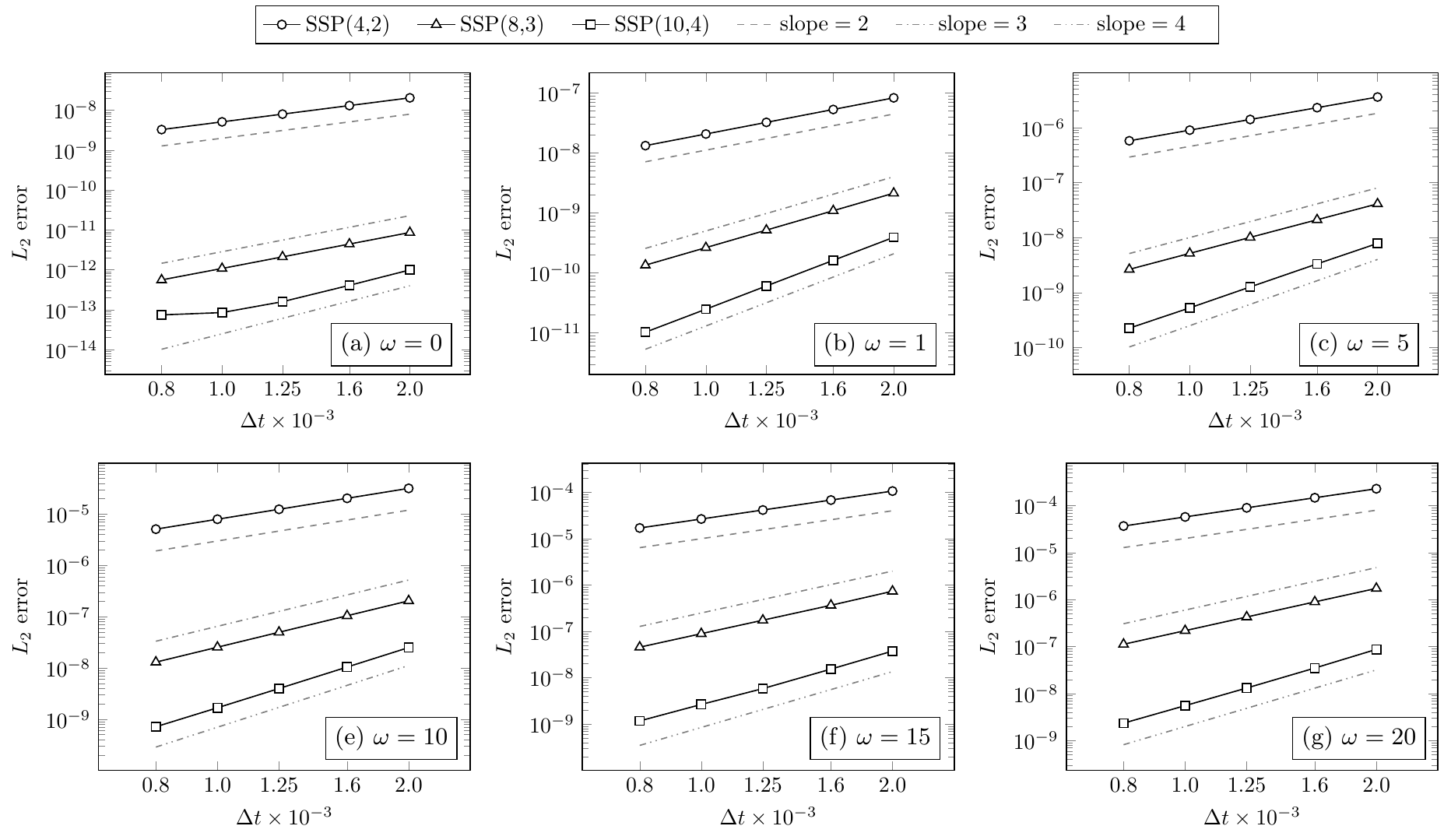}
\caption{$L_2$ errors of $\rho$ against time step size for an Euler vortex flow.}
\label{fig:euler_time}
\end{figure}

The tests are performed on the previous Euler vortex flow. Whenever it permits, we have varied the polynomial degree for each case to ensure that the spatial errors are negligibly small compared to the temporal errors. The results are reported in Fig. \ref{fig:euler_time} for different rotational speeds. For almost all the cases, the temporal errors decrease at the correct orders as the time step size decreases. For the $w=0$ case with SSP(10,4), the temporal errors become so small (in the order of $10^{-14}$) that they can no longer be easily separated from the spatial errors at small time step sizes, which results in the incorrect slope at small time step sizes. Nevertheless, at large time step sizes, the curve still shows the correct slope for this case. We also have tested the temporal accuracy on the Taylor-Couette flow, and similar results were obtained (for conciseness, the results are not included here).

\subsection{Conservation}
From the analysis in Sec. \ref{sec:conservation}, global conservation on a mesh with $N_{\text{elem}}$ elements, and $N_{\text{boun1}}$ and $N_{\text{boun2}}$ boundary cell faces that are mapped to the $\xi$ and the $\eta$ directions, respectively, can be expressed as
\begin{equation}
\sum_{e=1}^{N_{\text{elem}}}  \int_{0}^{1} \int_{0}^{1} \frac{\partial \widetilde{\mathbf{Q}}^e}{\partial t} \text{d}\xi\text{d}\eta +
\sum_{b_1=1}^{N_{\text{boun1}}} \int_{0}^{1} \text{sign}^{b_1} \cdot \widehat{\mathbf{F}}^{b_1} \text{d}\eta +
\sum_{b_2=1}^{N_{\text{boun2}}} \int_{0}^{1} \text{sign}^{b_2} \cdot \widehat{\mathbf{G}}^{b_2} \text{d}\xi = \mathbf{0},
\label{eq:conserve1}
\end{equation}
where $\text{sign}=1$ or $-1$, depending on which standard face the physical cell face is mapped to. Denoting the temporal discretization by {\small $\widetilde{\partial}/\widetilde{\partial} t$}, then the exact time derivative in (\ref{eq:conserve1}) can be expressed as
\begin{equation}
\frac{\partial \widetilde{\mathbf{Q}}^e}{\partial t} = \frac{\widetilde{\partial} \widetilde{\mathbf{Q}}^e}{\widetilde{\partial} t} + \boldsymbol{\epsilon}(\Delta t, \xi, \eta),
\end{equation}
where $\boldsymbol{\epsilon}(\Delta t, \xi, \eta)$ represents a small error that is a function of time-step size and space. Substitute the above relation into (\ref{eq:conserve1}),
\begin{equation}
\sum_{e=1}^{N_{\text{elem}}}  \int_{0}^{1} \int_{0}^{1} \frac{\widetilde{\partial} \widetilde{\mathbf{Q}}^e}{\widetilde{\partial} t} \text{d}\xi\text{d}\eta +
\sum_{b_1=1}^{N_{\text{boun1}}} \int_{0}^{1} \text{sign}^{b_1} \cdot \widehat{\mathbf{F}}^{b_1} \text{d}\eta +
\sum_{b_2=1}^{N_{\text{boun2}}} \int_{0}^{1} \text{sign}^{b_2} \cdot \widehat{\mathbf{G}}^{b_2} \text{d}\xi = \mathbf{E}(\Delta t),
\label{eq:conserve2}
\end{equation}
where
\begin{equation}
\mathbf{E}(\Delta t) = - \sum_{e=1}^{N_{\text{elem}}} \int_{0}^{1} \int_{0}^{1} \boldsymbol{\epsilon}^{e}(\Delta t, \xi, \eta) \text{d}\xi\text{d}\eta.
\end{equation}
When $\Delta t \to 0$ or when it is a steady-state problem, we will have $\mathbf{E}(\Delta t) \to 0$ if a numerical scheme is indeed conservative. And the four components of $\mathbf{E}(\Delta t)$ measure how well conservation is satisfied for mass, momentums and energy.

Based on the above analysis, we employ a flat-plate Couette flow to verify the conservation property of the present method. This flow is formed between two infinite flat plates separated by a distance of $H$, with the upper plate moving at a constant speed $U$ and the lower plate fixed. The steady state solution is
\begin{gather}
u = U \frac{y}{H}, \quad v = 0, \quad p = \text{constant}, \\[0.3em]
T = T_0 + \big(T_1 - T_0 \big) \frac{y}{H} + \frac{\mu U^2}{2 \kappa} \left(\frac{y}{H}-\frac{y^2}{H^2} \right),
\end{gather}
where $0\le y \le H$ is the vertical coordinate, $T_0$ and $T_1$ are temperatures on the lower and the upper plates.

The mesh used for this simulation is identical to that in Fig. \ref{fig:euler_msh}, but is scaled by a factor of $1/10$ in each direction (i.e., with $H=1$) to make the simulation setup easier. The upper plate speed is set to $U=1$. The temperatures are set to $T_0=T_1=1$. All the other parameters are chosen such that the flow has a Prandtl number $Pr=0.72$ and Reynolds number $Re=100$, and the flow on the upper plate has a Mach number $Ma=0.8$. The flow field is initialized using the exact solution. The boundary conditions are weakly imposed also using the exact solution. The SSP(5,4) temporal scheme with a time step size of $1.0\times 10^{-5}$ is used for all the simulations. And all the simulations are performed using double precision float numbers.

We have tested different rotational speeds and polynomials, with $\mathbf{E}(\Delta t)$  monitored until $Ut/H=20$. It was noticed that for all the cases, $\mathbf{E}(\Delta t)$ stays in the level of machine precision and does not grow during the simulation. For conciseness, we only report the values of $|\mathbf{E}(\Delta t)|$ for some of the tests in Tab. \ref{tab:conservation} at the end of each simulation. In the calculation, the numerical time derivative term in (\ref{eq:conserve2}) is replaced by the residual according to (\ref{eq:resid}), and the integrations are evaluated exactly using quadratures.
\begin{table}[!htb]
\def\arraystretch{1.1}
\setlength{\tabcolsep}{4mm}
\centering
\begin{tabular}{>{\raggedleft}m{4mm} >{\raggedright}m{7mm} | >{\centering}m{22mm} >{\centering}m{22mm} >{\centering}m{22mm} >{\raggedleft\arraybackslash}m{20mm}}
\hline
$\omega$ & $P$ & Mass & Momentum 1 & Momentum 2 & Energy~~~ \\
\hline
0  & \hphantom{1}4 & 1.300E-16 & 3.141E-16 & 4.935E-16 & 2.386E-16 \\
   & \hphantom{1}8 & 2.428E-16 & 1.103E-16 & 8.592E-16 & 3.013E-17 \\
   &            12 & 2.182E-16 & 5.174E-16 & 1.699E-15 & 7.907E-16 \\
\hline
5  & \hphantom{1}4 & 1.490E-17 & 1.837E-15 & 1.078E-16 & 2.797E-15 \\
   & \hphantom{1}8 & 4.948E-16 & 9.753E-16 & 9.121E-17 & 8.857E-16 \\
   &            12 & 7.329E-16 & 7.710E-16 & 1.142E-15 & 4.289E-15 \\
\hline
10 & \hphantom{1}4 & 2.308E-15 & 3.565E-15 & 1.919E-16 & 3.406E-15 \\
   & \hphantom{1}8 & 2.595E-16 & 2.406E-16 & 1.905E-16 & 1.256E-15 \\
   &            12 & 1.375E-15 & 1.132E-15 & 5.528E-16 & 8.048E-15 \\
\hline
20 & \hphantom{1}4 & 1.920E-15 & 1.615E-15 & 1.151E-15 & 3.046E-15 \\
   & \hphantom{1}8 & 1.606E-15 & 1.183E-16 & 3.572E-16 & 8.508E-15 \\
   &            12 & 5.550E-15 & 9.471E-16 & 1.081E-15 & 8.013E-15 \\
\hline
\end{tabular}
\caption{Conservation of mass, momentum and energy on a flat-plate Couette flow.}
\label{tab:conservation}
\end{table}

\subsection{Free-Stream Preservation}
In this test, the free stream flow is chosen such that it has a Mach number $Ma=0.3$. The same mesh as that for the conservation test are used for this test. Dirichlet boundary conditions using the constant flow are applied at all the outer boundaries. Different polynomials and rotational speeds are tested to investigate their effects on free-stream preservation. For all the cases, we employ SSP(10,4) with a time-step size of $1.0\times 10^{-3}$ as the time marching scheme. We measure the free-stream preservation by the normalized $L_2$ error of pressure (represents the average strength of numerically generated disturbances). As a comparison, we also test free-stream preservation on a dynamic conforming mesh that has almost exactly the same resolution as the nonconforming sliding mesh. The center of conforming mesh has a displacement of $\Delta y(t)=0.1 \sin t$ with respect to its initial position, with the mesh in the rest of the domain deforming using an algebraic function \cite{zhang-2016a}. For all the cases, we monitored the $L_2$ errors until the nondimensional time reaches $U_\infty t/H=20$. The results at the end of the simulations are summarized in Tab. \ref{tab:freestream1}.

For perfect free-stream preservation, the error should stay in the level of machine precision, and its magnitude should not change dramatically with polynomial degree. Table \ref{tab:freestream1} demonstrates that the FR method together with the GCL equations ensure very good free-stream preservation on the dynamic conforming mesh. The sliding mesh method, however, does not perfectly satisfy free-stream preservation. Instead, it approaches free-stream preserving in an exponential manner as the polynomial degree increases.
\begin{table}[!htb]
\def\arraystretch{1.3}
\setlength{\tabcolsep}{2mm}
\centering
\begin{tabular}{c| *{7}{c}}
\hline
\backslashbox{~~$\omega$}{$P$~~}     & 2 & 3 & 4 & 5 & 6 & 7 & 8  \\ \hline
$\ast$ & 6.213E-16 & 1.983E-15 & 3.869E-15 & 4.009E-15 & 3.566E-15 & 2.991E-15 & 2.955E-15 \\ \hline
0      & 7.844E-06 & 2.585E-07 & 8.252E-09 & 2.173E-10 & 4.676E-12 & 1.049E-13 & 3.332E-15 \\ \hline
1      & 7.609E-06 & 2.006E-07 & 7.495E-09 & 1.625E-10 & 3.985E-12 & 7.057E-14 & 3.213E-15 \\ \hline
5      & 7.642E-06 & 2.075E-07 & 7.393E-09 & 1.614E-10 & 4.006E-12 & 7.465E-14 & 3.518E-15 \\ \hline
10     & 7.559E-06 & 2.204E-07 & 7.283E-09 & 1.691E-10 & 4.036E-12 & 7.857E-14 & 3.364E-15 \\ \hline
15     & 7.524E-06 & 2.491E-07 & 7.806E-09 & 2.068E-10 & 4.455E-12 & 9.981E-14 & 3.658E-15 \\ \hline
20     & 7.817E-06 & 3.456E-07 & 8.873E-09 & 2.612E-10 & 5.123E-12 & 1.149E-13 & 3.880E-15 \\ \hline
\end{tabular}
\caption{Free stream preservation on a dynamic conforming mesh (marked by `$\ast$') and a sliding mesh (the rest).}
\label{tab:freestream1}
\end{table}

The reason for the present method not strictly satisfying free-stream preservation is evident. For a constant free stream flow, we are basically projecting the metric terms between cell faces and mortars. When a circular cell face is represented exactly, its metric terms are in the space of $\boldsymbol{\mathsf{P}}_\infty$. But the output of a projection is in the space of $\boldsymbol{\mathsf{P}}_{N-1}$. This means that there is always a truncation error of $\mathcal{O}(\xi^N)$ during the projection of metric terms, which therefore results in an exponential approximation of free-stream preservation.

The issue of using exact geometric expression for calculating metric terms has been discussed in \cite{kopriva-2006,kopriva-2009}, and the suggestion is to use polynomial approximation to satisfy free-stream preservation. This suggestion has been proven to work perfectly on curved conforming meshes. A recent work \cite{kopriva-2019} on curved nonconforming meshes in the special scenario of subdivision of parent element reveals that a necessary condition for free-stream preservation is that the metrics of a child face being computed from its parent face. But for general curved nonconforming meshes, such as the sliding meshes in this work, a child face (e.g., a mortar) does not have a unique parent (more specifically, a child face has two parent faces that are different polynomials), and thus the aforementioned necessary condition generally can not be satisfied. To the authors' knowledge, free-stream preservation on general curved nonconforming meshes still remains an open problem for all polynomial-based high-order methods. Note that for linear nonconforming meshes, such as that in \cite{qiu-2019}, free-stream preservation can always be easily satisfied. But linear mesh obviously introduces too large geometric error for the curved sliding meshes here (see Appendix \hyperref[sec:appendix_transfinite]{A}).

Nevertheless, Tab. \ref{tab:freestream1} shows that the free-stream-preservation errors from the sliding mesh quickly decreases to machine precision even at a moderate polynomial degree of 8. Since this test is done on a very coarse grid, we thus expect very minor free-stream preservation effects in real applications where the meshes are usually finer.

\subsection{Flow over a Rotating Square Cylinder}
\begin{figure}[!htb]
\centering
\includegraphics[width=0.40\textwidth]{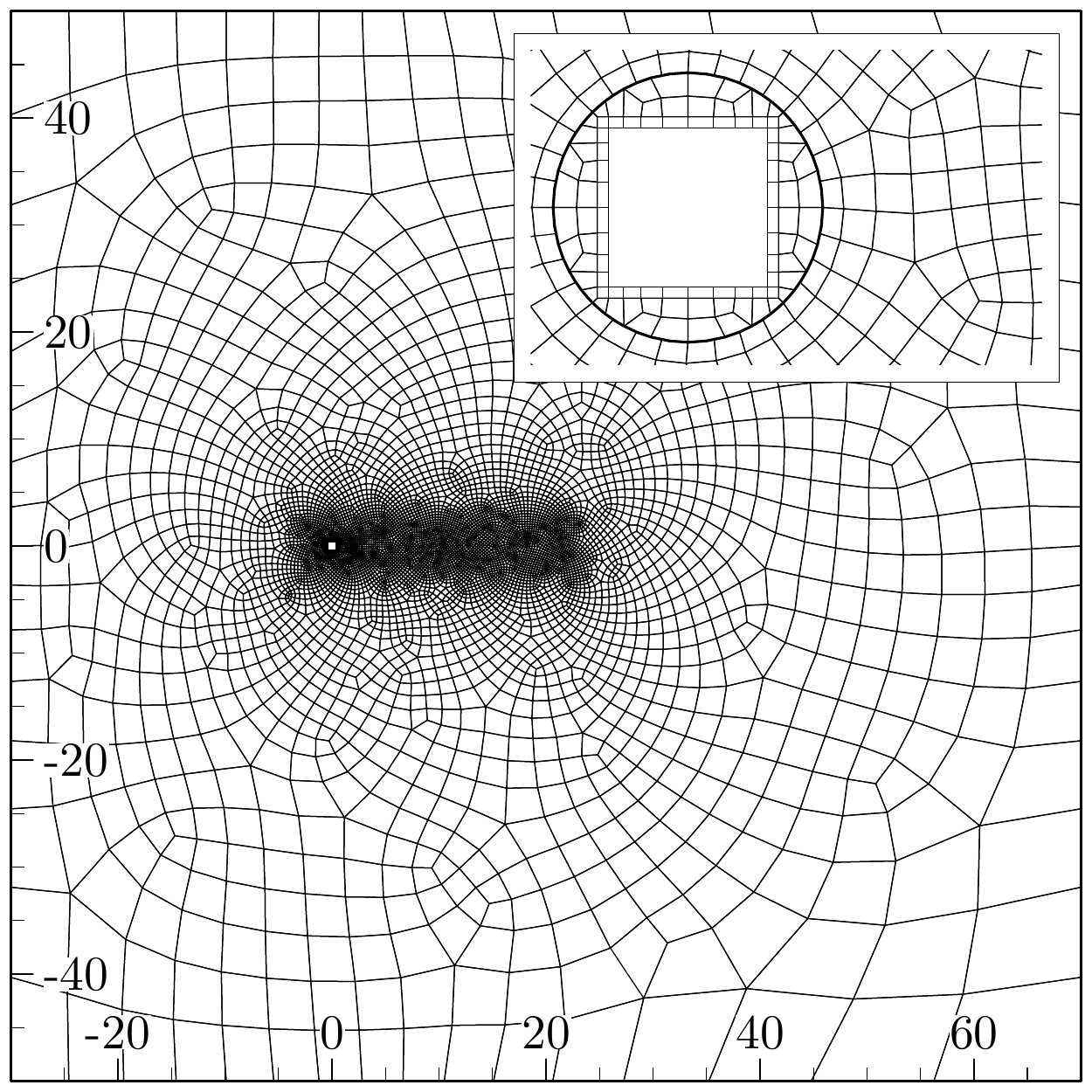} \qquad
\includegraphics[width=0.40\textwidth]{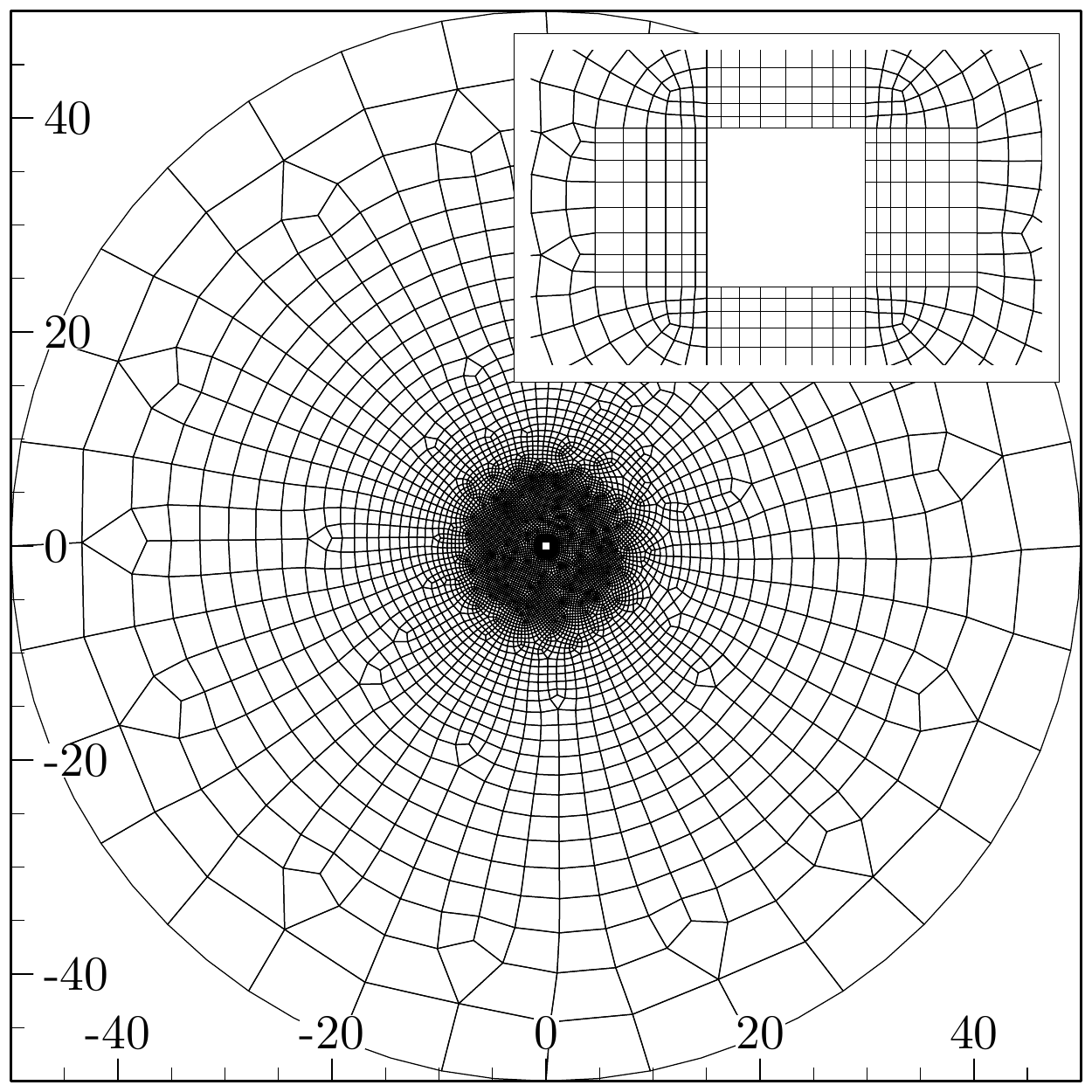}
\caption{Meshes for a rotating square cylinder: left, sliding mesh; right, rigid-rotating mesh.}
\label{fig:cylinder1msh}
\end{figure}

As shown in Fig. \ref{fig:cylinder1msh}, we employ two types of meshes for this simulation: a sliding mesh (on the left) and a rigid-rotating mesh (on the right). The square cylinder has a diameter (i.e., diagonal length) of $D=1$. The sliding mesh consists of 6,797 quadrilateral elements in an overall $100 D \times 100 D$ domain, with 95 of the elements inside a small rotating subdomain whose diameter is $1.2D$. The rigid-rotating mesh consists of 6,841 elements on a whole-piece rotating circular domain whose diameter is $100D$. Both meshes are refined in the vicinity of the cylinder to better resolve flow structures. With similar amounts of mesh elements, the sliding mesh obviously provides good resolution in a wider range of the wake region. The rigid-rotating mesh, on the other hand, always wastes a large amount of elements in unimportant regions, and this is unavoidable. Close views of the meshes are shown on the top right of each figure.

The freestream flow is chosen such that it has a Mach number $Ma=0.1$. The Reynolds number based on free-stream flow properties and the cylinder diameter is $Re_D=100$. For the sliding mesh, the inner subdomain rotates counterclockwise at a non-dimensional rotational speed of $\omega D/U_\infty=\pi/2$ (which corresponds to a nondimensional period of $T^*=4$); and for the rigid-rotating mesh, the whole domain rotates at this speed. No-slip adiabatic wall boundary condition is applied on the cylinder surface. Characteristic far-field boundary conditions are applied at all outer boundaries. The eighth-order (i.e., $P=7$) spatial scheme and the SSP(10,4) time marching scheme \cite{ketcheson-2008} with a nondimensional time-step size of $\Delta t U_\infty/D=2.0\times 10^{-4}$ are employed to run the simulations. Polynomial and mesh refinement studies have also been performed to have confirmed that the present setup well resolves the flow.

\begin{figure}[!htb]
\centering
\includegraphics[width=0.98\textwidth]{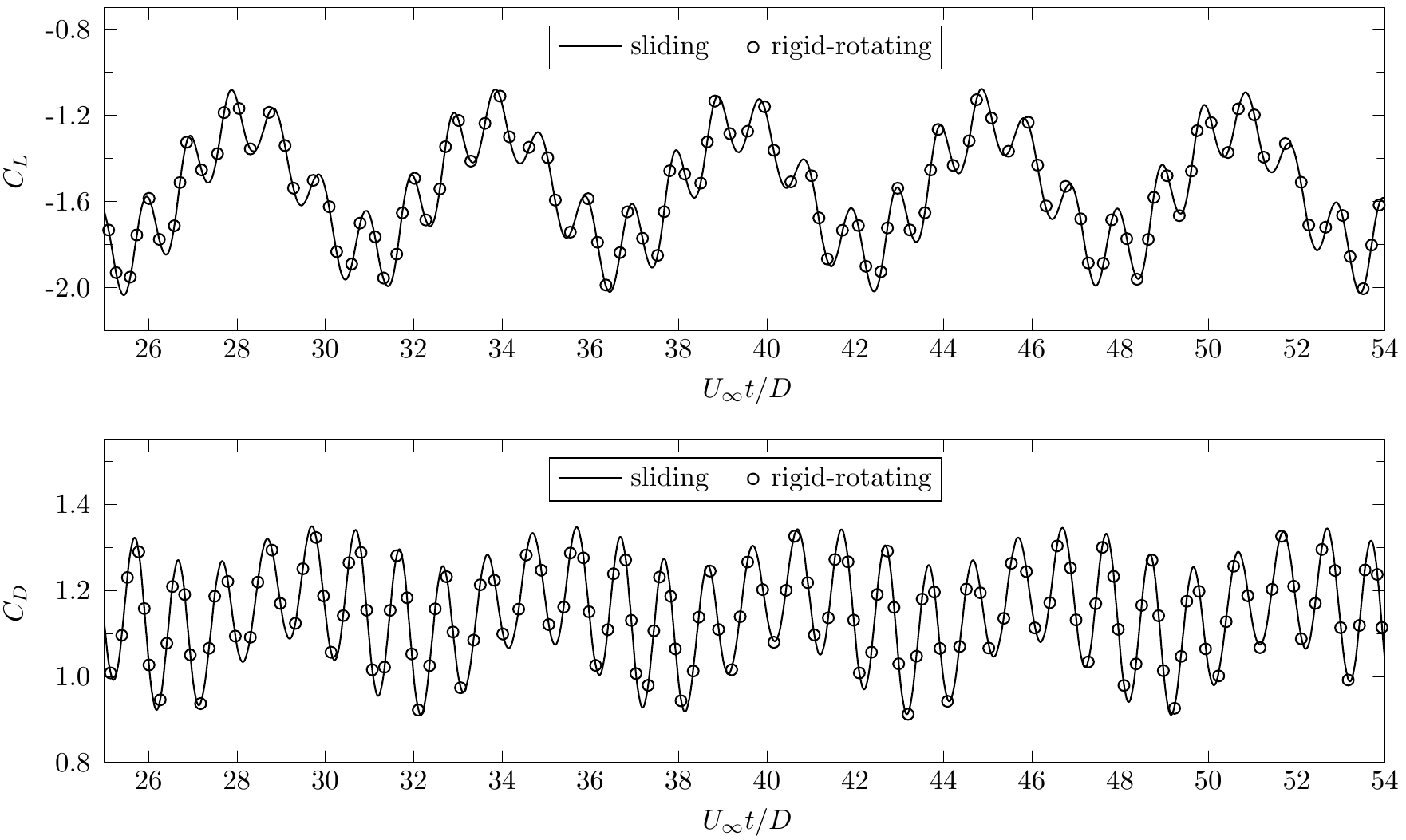}
\caption{Lift and drag coefficients of a rotating square cylinder using a sliding mesh and a rigid-rotating mesh.}
\label{fig:cylinder1clcd}
\end{figure}

The lift and drag coefficients of the cylinder are plotted in Fig. \ref{fig:cylinder1clcd}. The cylinder experiences a negative lift and positive drag all the time. These curves overall do not have a periodic pattern that is directly related to the cylinder's motion.
But the horizontal distance between two neighboring local peaks (or troughs) is approximately 1 (i.e., $T^*/4$), which is the period at which the cylinder disturbs the flow. Results from these two distinctly different approaches agree very well, which confirms the correctness of the present method.

\begin{figure}[!htb]
\centering
\includegraphics[scale=1.0]{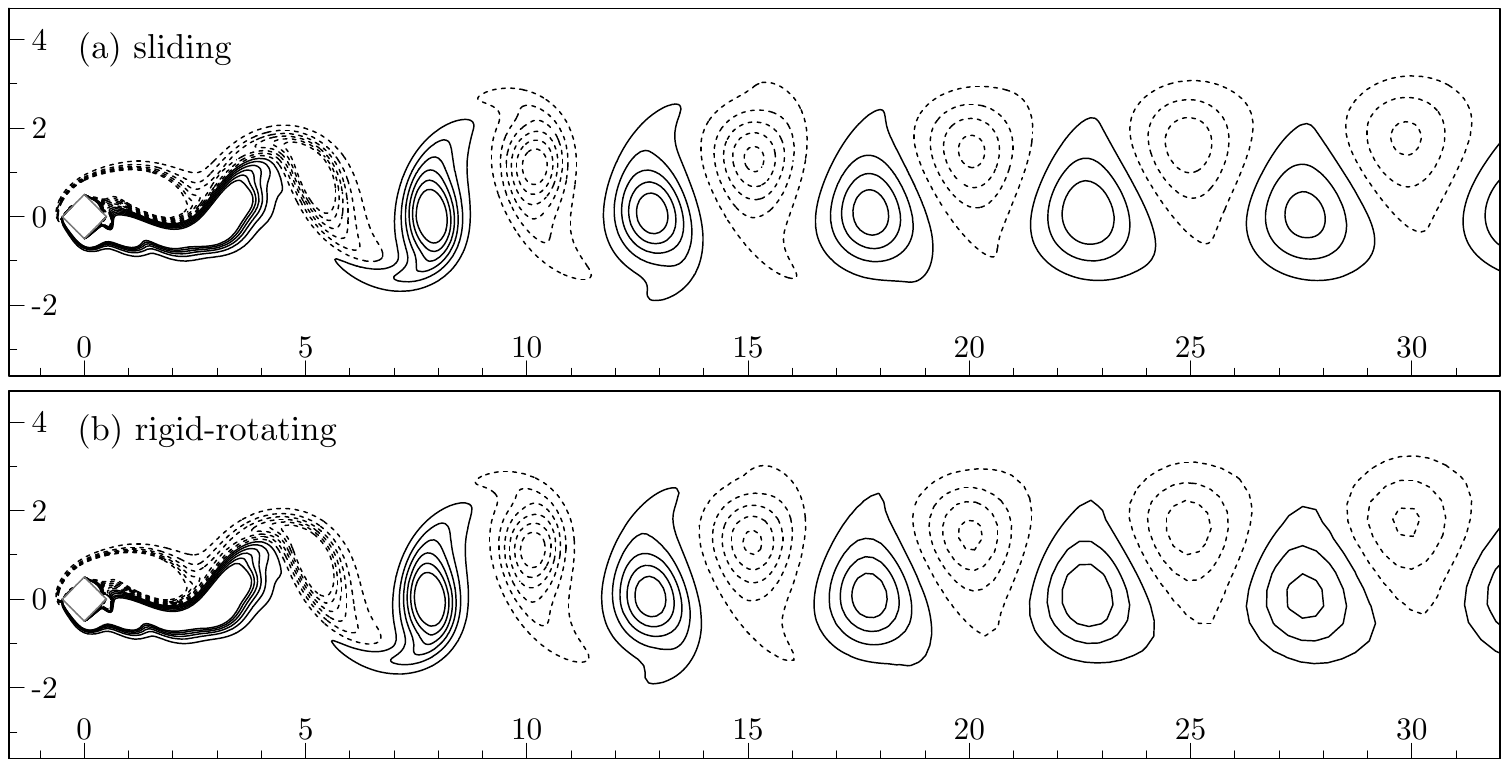}
\caption{Vorticity contours of flow over a rotating square cylinder (solid lines are positive, dashed lines are negative).}
\label{fig:cylinder1vort}
\end{figure}

We also compare the flow fields from the same time instant in Fig. \ref{fig:cylinder1vort} using vorticity contours. Alternative positive and negative vortices are observed in the wake region, and they are well captured even in the very far flow region. The two approaches basically produce identical flow fields as the vortices have almost exactly the same position and shapes. The far wake region from the sliding mesh is slightly better resolved than that from the rigid-rotating mesh. This is consistent with our previous observation that sliding mesh always provides better resolution than rigid-rotating mesh for a given number of elements.

\subsection{Flow over Multiple Rotating Square Cylinders}
\begin{figure}[!htb]
\centering
\includegraphics[scale=1.0]{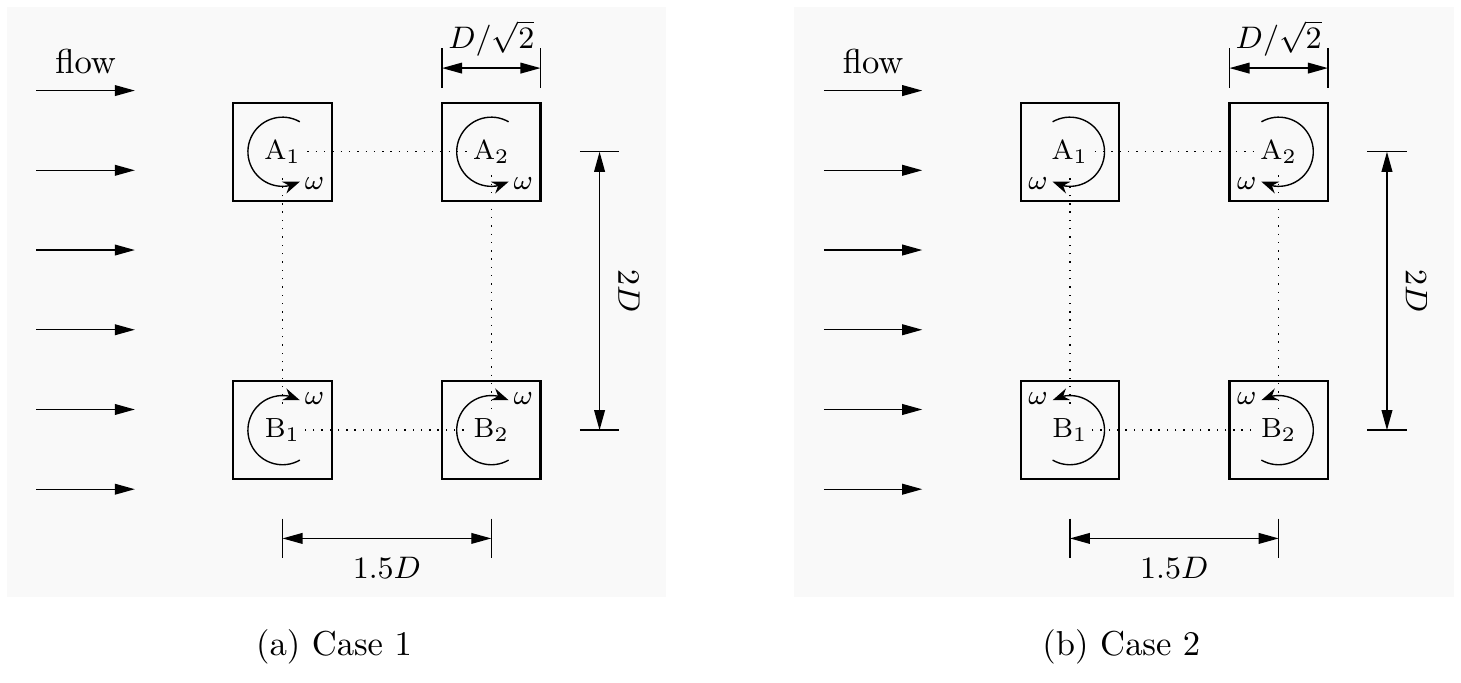}
\caption{Schematic of simulation setup for flow over multiple rotating square cylinders.}
\label{fig:cylinder2}
\end{figure}

In this test, we simulate the interactions of a uniform incoming flow and four rotating square cylinders to show the capability of the method for handling multiple rotating objects. A schematic of the simulation setup is shown in Fig. \ref{fig:cylinder2}. The four cylinders are separated by a horizontal distance of $1.5D$ and a vertical distance of $2D$, where $D$ is the diameter of the cylinder with the same definition as that in the previous test. Two cases are investigated, where the only difference is the rotational direction of each cylinder as marked in the schematic. For simplicity, we denote the four cylinders as A\textsubscript{1}, A\textsubscript{2}, B\textsubscript{1} and B\textsubscript{2}.

The overall mesh for this test has 7,409 elements, with its topology and distribution similar to those of the sliding mesh in the previous test. Each inner subdomain mesh is exactly the same as that in the previous test. For conciseness, views of the mesh are not repeated here. All the other parameters for this test, such as the freestream Mach number, the Reynolds number, the rotational speed, the boundary conditions, etc., are exactly the same as those for the previous test.  The simulations are performed using the eighth-order scheme with SSP(10,4) and a nondimensional time step size of $1.25\times 10^{-4}$. Each simulation is continued for a nondimensional time of 500. In what follows, we briefly discuss the results for each case.

\subsubsection{Case 1}
\begin{figure}[!htb]
\centering
\includegraphics[width=0.98\textwidth]{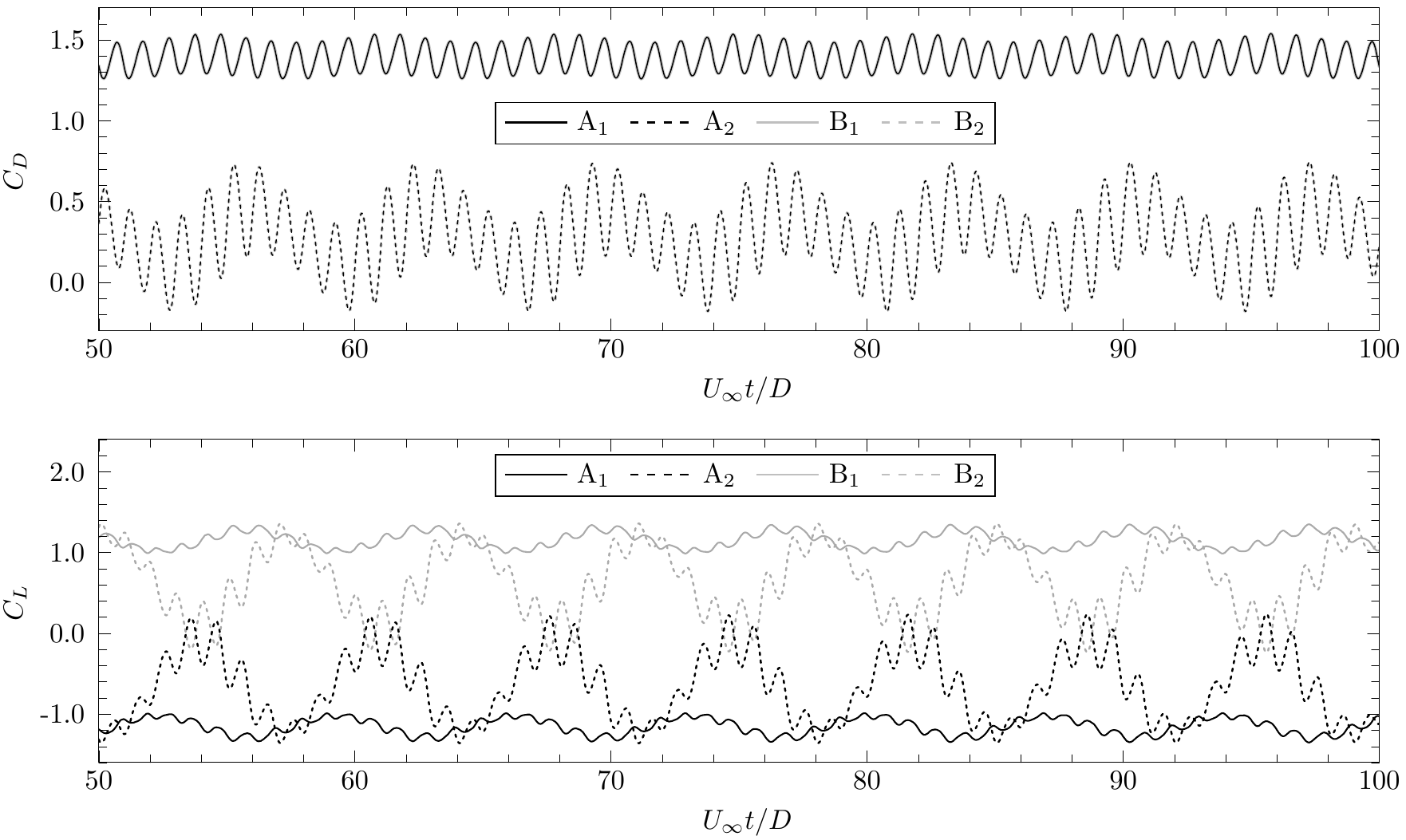}
\caption{Drag and lift coefficients for flow over multiple rotating square cylinders (case 1).}
\label{fig:cylinder2_clcd1}
\end{figure}

Figure \ref{fig:cylinder2_clcd1} shows the drag and lift coefficients of the cylinders from $U_\infty t/D=50$ to $100$. It is seen that cylinders from the same column experience the same drag in this time frame. In fact, we have observed from the simulation that this relation holds until $U_\infty t/D$ approaches $130$ when the curves start differing (not shown here, reasons explained later). The upstream cylinders experience much larger drag than the downstream ones do. The downstream ones, however, experience much larger drag oscillations, which is possibly due to vortex shedding. On the other hand, due to the Magnus effects and the opposite rotational directions, the lifts on A\textsubscript{1} and B\textsubscript{1} (also A\textsubscript{2} and B\textsubscript{2}) have the same magnitude but opposite signs. Again, the upstream cylinders experience smaller lift oscillations than the downstream cylinders do.

\begin{figure}[!htb]
\centering
\includegraphics[scale=1.0]{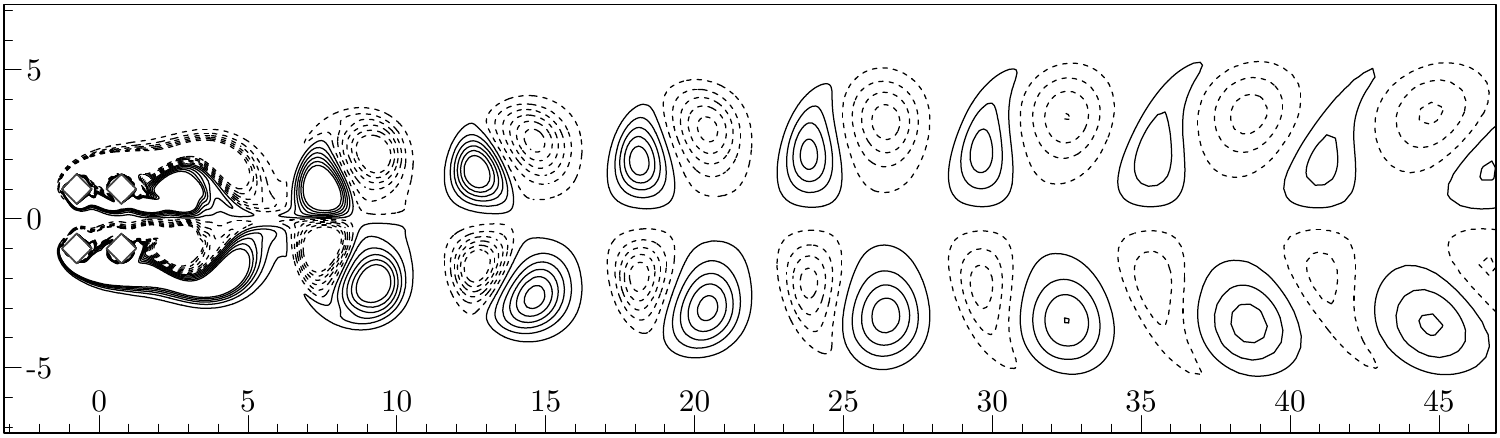}
\caption{Vorticity contours of flow over multiple rotating square cylinders (case 1).}
\label{fig:cylinder2_vort1}
\end{figure}

A snapshot of the flow field visualized by vorticity contours at $U_\infty t/D=80.5$ is shown in Fig. \ref{fig:cylinder2_vort1}. Well organized positive-negative vortex pairs are in the wake region and extend all the way to the very far flow field. This patten is perfectly mirror-symmetric about the horizontal center line of the setup. In fact, this is a very unstable and unsustainable system. Any tiny numerical disturbance could cause the flow to lose this symmetry. Using the present method, this symmetry is maintained for a surprisingly long time ($U_\infty t/D=130$) before breaking up. This case evidently demonstrates that the present method introduces very little numerical disturbance to a simulation. After the breakup of the symmetry, the flow becomes less interesting, and that is why we only focus on the flow before that point.

\subsubsection{Case 2}
\begin{figure}[!htb]
\centering
\includegraphics[width=0.98\textwidth]{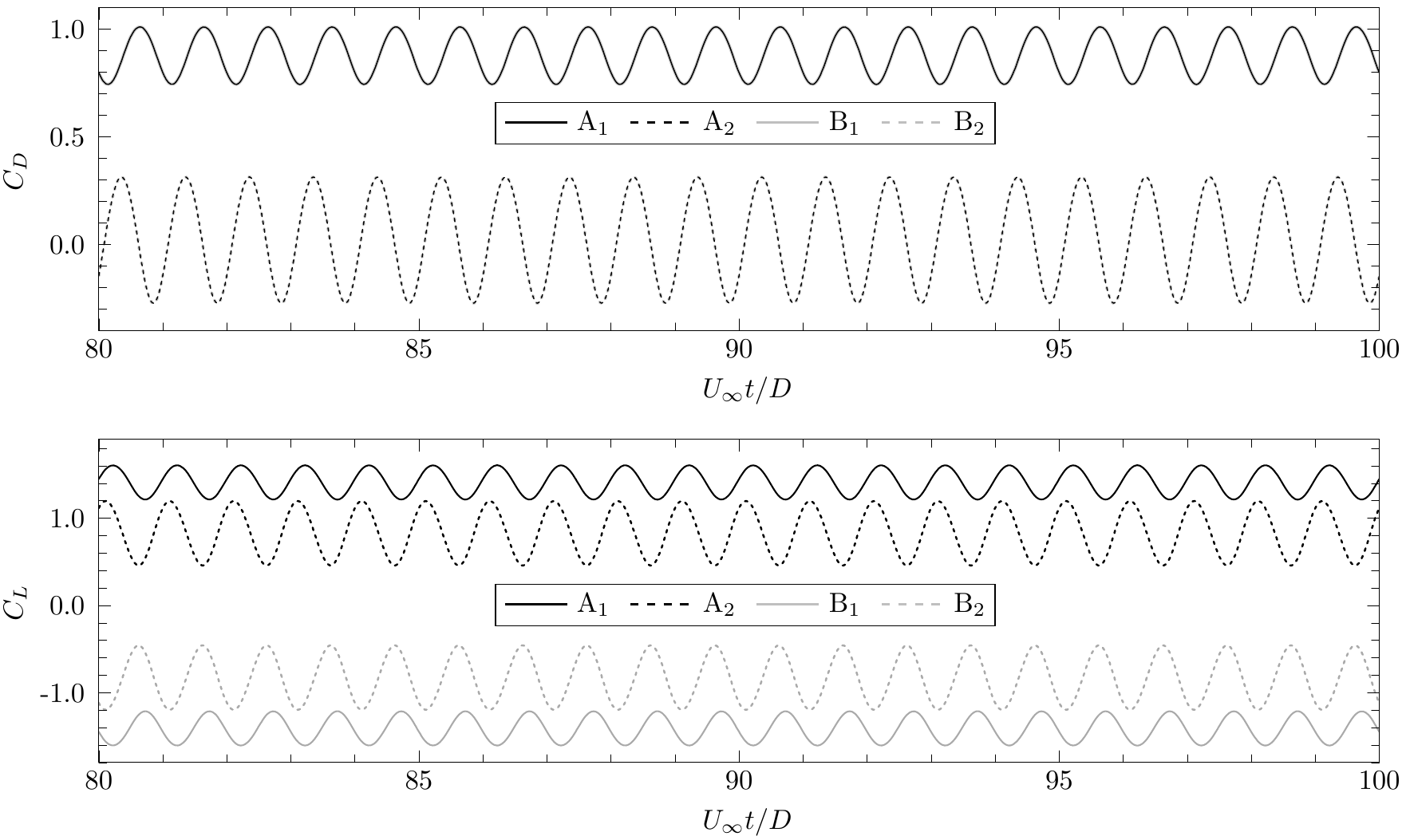}
\caption{Drag and lift coefficients for flow over multiple rotating square cylinders (case 2).}
\label{fig:cylinder2_clcd2}
\end{figure}

The drag and lift coefficients for this case are plotted in Fig. \ref{fig:cylinder2_clcd2}. The previous conclusions for Case 1 still hold here. Except that the current curves show a rather mono frequency that corresponds to a period of $T^*/4$, where $T^*$ is the nondimensional rotational period of the cylinders. This mono frequency indicates fewer vortical structures (or noises) in the flow field. Compared with Case 1, we also notice dramatic drag reduction on all the cylinders. On average, the rear cylinders even do not experience much drag at all. Meanwhile, dramatic lift enhancement is observed on all the cylinders.

\begin{figure}[!htb]
\centering
\includegraphics[scale=1.0]{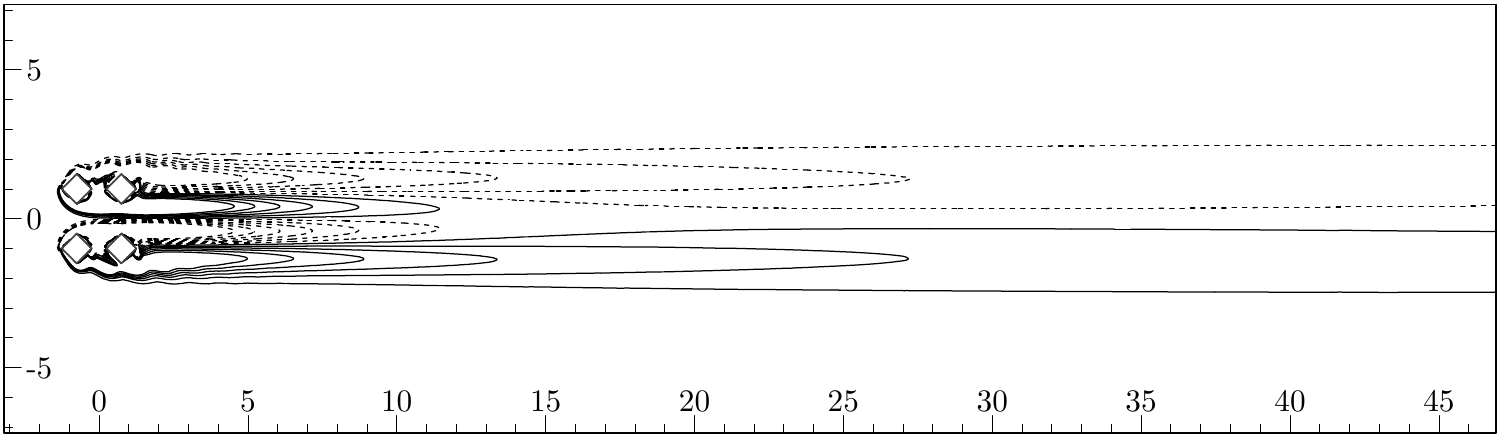}
\caption{Vorticity contours of flow over multiple rotating square cylinders (case 2).}
\label{fig:cylinder2_vort2}
\end{figure}

Figure \ref{fig:cylinder2_vort2} shows the instantaneous flow field at $U_\infty t/D=250.5$. Compared with Case 1, vortex shedding has been completely suppressed in this case, which is consistent with our observation form the force curves. This flow remains very stable and is almost steady (except in the very vicinity of the cylinders) throughout the simulation (i.e., even at $U_\infty t/D=500$). The reason for this steadiness is obvious. The rotation motion of the cylinders in this case (as shown in Fig. \ref{fig:cylinder2}(b)) decelerates the flow (also increases pressure) in the gap between the two rows. This reduction of momentum and increase of pressure prevent vortex shedding from the cylinders and thus make the flow almost steady.

\section{Summary}
\label{sec:conclusion}
We have presented a high-order sliding mesh method and performed detailed studies of its properties. This method employs transfinite mortar elements for communication between the two sides of a nonconforming sliding interface. Solutions and fluxes are projected back and forth in a least-squares fashion between mortars and cell faces to retain many of the favorable properties of the original conforming-mesh method. The use of transfinite mortar elements completely eliminates geometric errors on a circular sliding interface, and thus excludes an important source of errors from the scheme. Numerical tests have verified that this method retains the high-order accuracy of the FR method on both inviscid and viscous flows as the errors decrease exponentially with the increase of polynomial degrees. Numerical tests also have revealed that this method retains the order of accuracy of a time marching scheme, which allows the minimization of rotational speed effects by equipping the method with high-order time marching schemes. Through numerical analysis and verifications, we have confirmed the conservation property of the method, even under dramatically different rotational speeds. The satisfaction of outflow condition is also proved and verified, which ensures that, using the present method, a nonconforming sliding interface does not affect the characteristics of waves in a flow field. This method is, however, not strictly free-stream preservative. Instead, it approaches free-stream preservation in an exponential way. This is in fact not specific to the present method, as free-stream preservation on general curved nonconforming mesh still remains an open problem in the high-order method community. Nevertheless, we have shown through numerical tests that the free-stream-preservation errors are negligibly small even on a very coarse mesh using polynomials of moderate degrees. The comparison study of sliding mesh and rigid-rotating mesh on a rotating square cylinder has shown that the two approaches generate almost identical results. But a sliding mesh generally provides better wake resolution than a rigid-rotating mesh for a given number of overall mesh elements. The simulation of a matrix of rotating square cylinders has demonstrated the method's capability of simultaneously handling multiple objects. As a concluding remark, we must also emphasize that this method is not specifically designed for the FR method, it is readily applicable to many other high-order methods as well. This method is also easily extensible to deal with three-dimensional problems.

\section*{Acknowledgment}
This work was supported by a grant from the Office of Naval Research, administrated by Dr. Ki-Han Kim. The authors would like to express our acknowledgments for the support.

\newpage
\appendix

\section{Comparison of iso-parametric mapping and transfinite mapping}
\label{sec:appendix_transfinite}
\begin{figure}[!htb]
\centering
\includegraphics[scale=0.92]{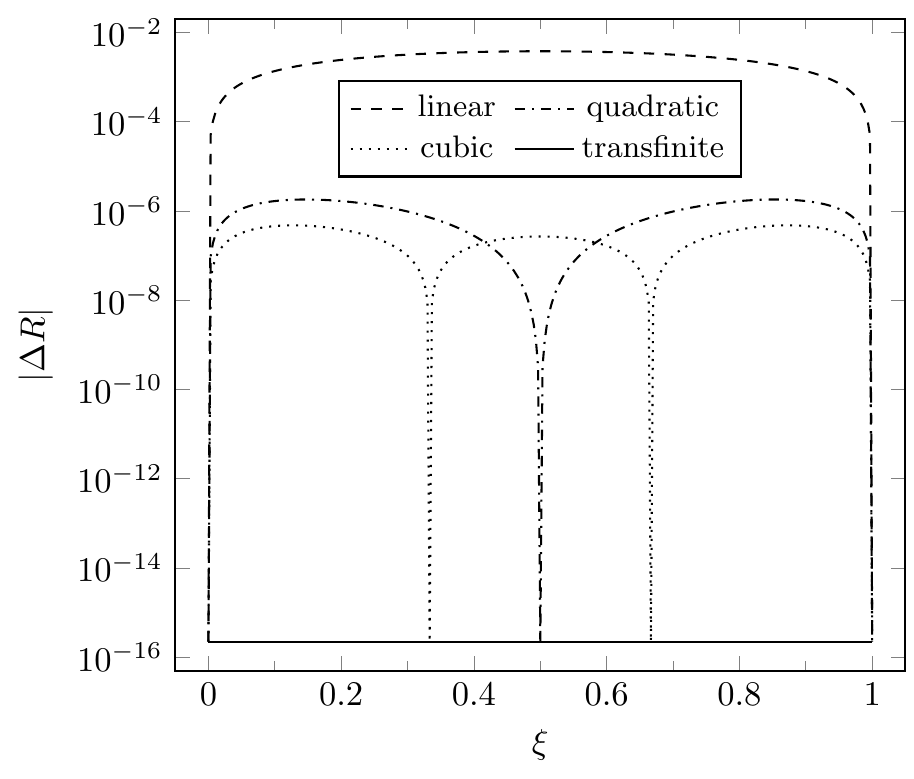} \\[-0.7em]
\caption{Comparison of iso-parametric mapping and transfinite mapping for approximating a circular arc.}
\label{fig:transfinite1d}
\end{figure}

Transfinite mapping can minimizes geometric errors, whereas iso-parametric mapping always carries a truncation error but monolithically approaches transfinite mapping as the order increases. We compare these two types of mappings through a 1D and a 2D examples. The 1D example is a circular arc with  $\theta_1=0^\circ$, $\theta_2=10^\circ$, $R=1$, and $\mathbf{x}_c=(0,0)$. We calculate the radius difference $\Delta R$ between the exact radius and the numerical ones that are calculated using the coordinates from the mappings. Figure \ref{fig:transfinite1d} shows that the iso-parametric mapping represents the circular arc exactly only at the nodal points. For non-nodal points, the error overall decreases as the order increases. On the other hand, the transfinite mapping always represents the circular arc exactly with the errors alway in the level of machine precision.

\begin{figure}[!htb]
\centering
\includegraphics[width=\textwidth]{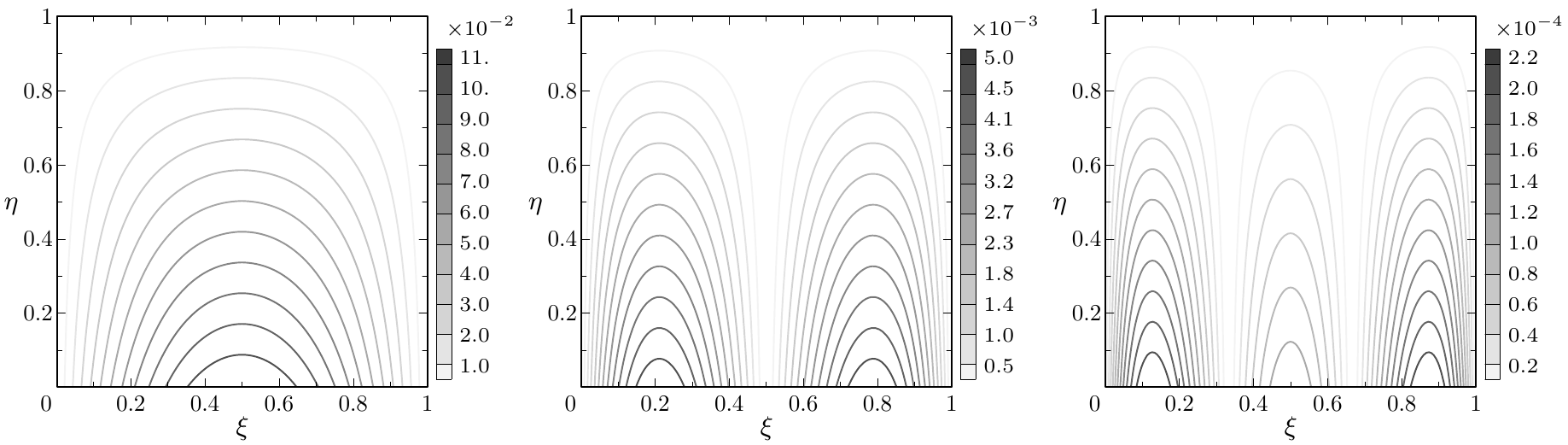} \\[-0.7em]
\caption{Difference between iso-parametric and transfinite mappings: left, $\Vert \mathbf{x}_\mathrm{L} - \mathbf{x}_\mathrm{T} \Vert$; middle, $\Vert \mathbf{x}_\mathrm{Q} - \mathbf{x}_\mathrm{T} \Vert$; right, $\Vert \mathbf{x}_\mathrm{C} - \mathbf{x}_\mathrm{T} \Vert$.}
\label{fig:transfinite2d}
\end{figure}

The 2D geometry is similar to that in Fig. \ref{fig:transfinite}, i.e., with three straight edges and one circular-arc edge. More specifically, we have chosen: $\mathbf{x}_1 = (R\cos\theta_1 + x_c, R\sin\theta_1 + y_c)$, $\mathbf{x}_2 = (R\cos\theta_2 + x_c, R\sin\theta_2 + y_c)$, $\mathbf{x}_3 = ((R/2)\tan\theta_2, R/2)$, $\mathbf{x}_4 = ((R/2)\tan\theta_1, R/2)$, where $R=2$, $\theta_1=-110^\circ$, $\theta_2=-70^\circ$, and $(x_c,y_c) = (0,R)$. Denote the coordinates from linear, quadratic, cubic, and transfinite mappings as $\mathbf{x}_\text{L}$, $\mathbf{x}_\text{Q}$, $\mathbf{x}_\text{C}$, and $\mathbf{x}_\text{T}$, respectively. There is no direct way to measure which of these coordinates are more accurate. Instead, we take the transfinite mapping as a reference, and plot the contours of the $L_2$ norm of coordinate difference (i.e., $\Vert \mathbf{x}_{*} - \mathbf{x}_{\text{T}} \Vert$, where `$*$' represents `L', or `Q', or `C') in Fig. \ref{fig:transfinite2d}. The minimal differences are seen around the nodal points and along the straight edges. The reason is that these geometric components are represented in the same way in all the mappings including the transfinite mapping. Overall, the $L_2$ norm decreases as the order of the iso-parametric mapping increases, i.e., the iso-parametric mapping approaches the transfinite mapping as the order increases.

\newpage
\biboptions{sort&compress}
\bibliographystyle{elsarticle-num-names}
\bibliography{references}



\end{document}